\documentclass[pdflatex,sn-chicago]{sn-jnl}

\usepackage{geometry}
\usepackage[utf8]{inputenc} 
\usepackage[T1]{fontenc}     
\usepackage[english]{babel}
\usepackage{dsfont}
\usepackage{cases}
\usepackage{makecell}
\usepackage{lineno}

\theoremstyle{thmstyleone}%
\newtheorem{theorem}{Theorem}
\newtheorem{proposition}[theorem]{Proposition}%

\theoremstyle{thmstyletwo}%
\newtheorem{remark}{Remark}%

\theoremstyle{thmstylethree}%

\def\thefootnote{\alph{footnote}}

\begin{document}

\title[Bayesian high-dimensional covariate selection in NLMEM]{Bayesian high-dimensional covariate selection in non-linear mixed-effects models using the SAEM algorithm}


\author*[1,2]{\fnm{Marion} \sur{Naveau}}\email{marion.naveau@inrae.fr}

\author[2]{\fnm{Guillaume} \sur{Kon Kam King}}

\author[3]{\fnm{Renaud} \sur{Rincent}}

\author[1]{\fnm{Laure} \sur{Sansonnet}}

\author[2]{\fnm{Maud} \sur{Delattre}}

\affil*[1]{\orgdiv{Université Paris-Saclay, AgroParisTech, INRAE}, \orgname{UMR MIA Paris-Saclay}, \orgaddress{\postcode{91120}, \city{Palaiseau}, \country{France}}}

\affil[2]{\orgdiv{Université Paris-Saclay, INRAE}, \orgname{MaIAGE}, \orgaddress{\postcode{78350}, \city{Jouy-en-Josas}, \country{France}}}

\affil[3]{\orgdiv{Université Paris-Saclay, INRAE, CNRS, AgroParisTech}, \orgname{GQE - Le Moulon}, \orgaddress{\postcode{91190}, \city{Gif-sur-Yvette}, \country{France}}}

\hypersetup{
pdftitle={Bayesian high-dimensional covariate selection in NLMEM using the SAEM algorithm},
pdfauthor={Marion ~Naveau, Guillaume ~Kon Kam King, Renaud ~Rincent, Laure ~Sansonnet, Maud ~Delattre},
pdfkeywords={High-dimension, Non-linear mixed-effects models, SAEM algorithm, Spike-and-slab prior, Variable selection},
}

\abstract{High-dimensional variable selection, with many more covariates than observations, is widely documented in standard regression models, but there are still few tools to address it in non-linear mixed-effects models where data are collected repeatedly on several individuals. In this work, variable selection is approached from a Bayesian perspective and a selection procedure is proposed, combining the use of a spike-and-slab prior and the Stochastic Approximation version of the Expectation Maximisation (SAEM) algorithm. Similarly to Lasso regression, the set of relevant covariates is selected by exploring a grid of values for the penalisation parameter. The SAEM approach is much faster than a classical MCMC (Markov chain Monte Carlo) algorithm and our method shows very good selection performances on simulated data. Its flexibility is demonstrated by implementing it for a variety of nonlinear mixed effects models. The usefulness of the proposed method is illustrated on a problem of genetic markers identification, relevant for genomic-assisted selection in plant breeding.}


\pacs[Acknowledgement]{
\rightskip=0pt\leftskip=0pt
The authors thank the Stat4Plant project ANR-20-CE45-0012 for funding. The authors are grateful to the INRAE MIGALE bioinformatics facility (MIGALE, INRAE, 2020. Migale bioinformatics Facility, doi: 10.15454/1.5572390655343293E12) for providing computing and storage resources. All our experiments have been done on this platform with the OPENBLAS library. The authors are grateful to Emmanuel Heumez (INRAE Estrées-Mons, UE GCIE) for managing the trial in Estrées-Mons, and to Jacques Le Gouis for managing the PIA BreedWheat. The data-set was obtained thanks to the support of the PIA (Investment for the Future Program) Breedwheat (ANR-10-BTBR-03).
}


\maketitle

\section{Introduction}

Mixed-effects models have been introduced to analyse observations collected repeatedly on several individuals in a population of interest \citep{lavielle, pinheiro}. This type of data is particularly common in the fields of pharmacokinetics or when modelling biological growth for example, where data is customarily analysed using a nonlinear model whose coefficients often have a biological interpretation. In this case, the intrinsic variability of the data captured by the model parameters is then attributable to different sources (intra-individual, inter-individual, and residual) whose consideration is essential to characterise without bias the biological mechanisms behind the observations. Mixed-effects models allow the study of the responses of individuals with the same overall behaviour but with individual variations characterised by random individual parameters that are not observed. Thus, mixed-effects models are latent variable models. Parameter inference is therefore difficult because the likelihood and classical estimators do not have an explicit form. A widely used solution is to use an EM (Expectation-Maximisation) algorithm, or any variant, to compute the maximum likelihood estimator or the maximum \textit{a posteriori} estimator in a Bayesian framework \citep{dempster}. 

Moreover, the description of inter-individual variability may involve a number of covariates much larger than the number of individuals. In this high-dimensional context, it is often desirable to be able to focus on the few most relevant covariates through a variable selection procedure. However, in mixed-effects models, identifying the influential covariates is difficult, as the selection concerns latent variables in the model. Recent years have seen the emergence of varied contributions on high-dimensional covariate selection in mixed-effects models. 
The proposed tools are very different according to whether the regression function is linear or non-linear with respect to the individual parameters.
More precisely, the linear case allows the development of criteria whose calculation and/or theoretical study involve explicit quantities, which is rarely true when the model is non-linear. 
In linear mixed-effects models, many rely on the use of regularised methods (see \cite{schelldorfer2011} and \cite{fan} for example) and most of them include theoretical consistency results that guarantee the good properties of the proposed methods. 
In contrast, in the more general framework of non-linear mixed-effects models (NLMEM), there are few results and the only published works concern computational aspects. \cite{bertrand} compare a stepwise approach using an empirical Bayes estimate and penalised regression approaches like Ridge, Lasso and HyperLasso penalties to select covariates in pharmacokinetics studies, but without calibrating the penalty parameters. \cite{ollier} proposes a proximal gradient algorithm for computing a Lasso estimator in order to perform joint selection of covariates and correlation parameters between random effects, this method being more suitable in low rather than large dimensional covariate settings. The codes available for the above-mentioned methods are not generic and tuning them is difficult. So the easiest practical solution for a practitioner faced with a high-dimensional problem is to adopt a two-stage approach, \textit{i.e.} to fit independent nonlinear models to each individual and perform variable selection using the estimated parameters, losing the uncertainty on these estimates and the beneficial shrinkage property of mixed-effect models.

Bayesian approaches to variable selection have not received much attention in the NLMEM context. The focus for their development has been classical statistical models like linear regression or generalised linear model, for which Bayesian variable selection has been intensively developed in recent years. 
These methods encourage sparsity in the regression vector by using a variety of priors (see for example \cite{tadesse2021handbook} and the references therein), which may have better properties than the double-exponential prior associated with the Lasso penalty. 
Very recently, \cite{lee_bayesian_2022} proposed an overview of the formulation, interpretation and implementation of Bayesian non-linear mixed-effects models. 
In particular, he discussed Bayesian inference methods, priors options, and model selection methods in this context.
However, these Bayesian approaches are based on Markov Chain Monte-Carlo (MCMC) methods which seldom scale well enough to be usable for high-dimensional variable selection. 
The main objective of this paper is to propose a fast Bayesian spike-and-slab approach that can be used to identify the relevant covariates in a non-linear mixed-effects model in a high-dimensional context. 
More precisely, we extend the EMVS (EM Variable Selection) approach of \cite{rockova} to the NLMEM setting.
Like EMVS, the proposed approach involves two major steps. The first step is, for different values of the spike hyperparameter, to select a local version of the median probability model \citep{barbieri2004optimal} using
the Stochastic Approximation version of the EM algorithm (SAEM, see \cite{delyon} and \cite{kuhn}). 
The second step consists in selecting the "best" model among those kept after the first step, using an extension of the BIC criterion \citep{chen}.
An important difference with \cite{rockova} is that our approach is applied to NLMEM and not to classical linear regression models. Due to the model non-linearity and to the latent nature of the model random effects, the central so-called $Q$-quantity of the EM algorithm often does not have a closed-form expression and posterior distributions are difficult to compute. To overcome these issues, we propose an inference method using the SAEM algorithm rather than simply the EM algorithm as in \cite{rockova}. Another important difference is that optimal model selection among the sub-models obtained in the first step does not require the calculation of the marginal posterior of the models for a spike parameter being equal to 0, as in \cite{rockova}, but only of the log-likelihood of the NLMEM taken at the maximum likelihood estimator.

This paper is organised as follows. Section~\ref{stat_model} describes the non-linear mixed-effects model to introduce the notation, summarises the main objective of our procedure, and defines and motivates the hierarchical prior formulations. Section~\ref{inference} details the key tools for our approach: the SAEM algorithm, to compute the maximum \textit{a posteriori} estimator of the model parameters and a thresholding rule to select a local version of the median probability model and put some coefficients of the regression vector to zero. Next, Section~\ref{Model_select} describes the variable selection procedure. Section~\ref{simul} evaluates the selection performance of our method through an intensive simulation study and presents comparisons with existing methods. To illustrate the high-dimensional variable selection method proposed in this paper, Section~\ref{Section_data_reel} presents an application on a real data-set composed of European elite winter wheat varieties. Finally, Section~\ref{conclusion} concludes with a summary discussion and prospects for future research. More details about algorithms and real data are postponed in the appendices. 

\section{Model description}
\label{stat_model}

\subsection{Nonlinear mixed-effects model and notations}
\label{notations}

The formalism used in this paper is that of \cite{lavielle} and \cite{pinheiro}. Let $n$ be the number of individuals and $n_i$ the number of observations for individual~$i$. Let $n_{tot}=\sum_{i=1}^{n} n_i$ the total number of observations. Let $X^\top$ denote the transpose of a vector or matrix $X$. Consider the following non-linear mixed-effects model: for all $1\leq i \leq n$ and $1\leq j \leq n_i$,
\begin{subnumcases}{\label{model} }
    y_{ij} = g(\varphi_i, t_{ij}) + \varepsilon_{ij}, & $\varepsilon_{ij} \overset{\text{i.i.d.}}{\sim} \mathcal{N}(0,\sigma^2)$, \label{niv1}\\
    \varphi_i = \mu + \boldsymbol{\beta}^\top V_i+ \xi_i, & $\xi_i \overset{\text{i.i.d.}}{\sim} \mathcal{N}_q(0,\Gamma)$. \label{niv2}
\end{subnumcases} 
This model is described in two levels. 
First, at the individual level, Equation~\eqref{niv1} describes the intra-individual variability, where the observations $y_{ij}$ in $\mathbb{R}$ represent the response of individual $i$ at time $t_{ij}$. It is assumed that all individuals follow the same known functional form $g$ which depends non-linearly on an individual parameter $\varphi_i$ which is $q$-dimensional. Thus, this function governs intra-individual behaviour. The value of $q$ is closely linked to the choice of the mechanistic function $g$. It is therefore known, and usually small. The variance $\sigma^2>0$ of the Gaussian measurement noise is assumed unknown. Then, at the population level, Equation~\eqref{niv2} describes the inter-individual variability. For all $i \in \{1,\dots, n\}$, the $q$ individual parameters $\varphi_{i}=(\varphi_{im})_{1\leq m\leq q} \in \mathbb{R}^q$ are modelled as a multivariate Gaussian random variable whose mean is specified as the sum of an intercept $\mu$ in $\mathbb{R}^q$ and a linear combination of known covariates measured on individual $i$ and contained in the vector $V_{i}=(V_{i1}, \dots, V_{ip})^\top \in \mathbb{R}^p$. The term "covariates" refers to explanatory variables which may be relevant to explain inter-individual variability. This term is used to distinguish them from other explanatory variables such as the time variable for example. The number of covariates is denoted by $p$, $\boldsymbol{\beta}=(\beta_{\ell m})_{1\leq \ell \leq p ; 1\leq m \leq q} \in \mathcal{M}_{p \times q}$ is an unknown fixed effects matrix, and the inter-individual variance-covariance matrix $\Gamma$ is assumed unknown. Thus, the inter-individual variability is separated into two parts: on the one hand, $\boldsymbol{\beta}$ models the variability that can be explained by the covariates $V_i$, and on the other hand, $\xi_i$ represents the part of variation that is not explained by the measured covariates. In the following,
$y_i=(y_{ij})_{1 \leq j \leq n_i}$, $\boldsymbol{y}=(y_i)_{1 \leq i \leq n}$, $\boldsymbol{\varphi}=(\varphi_{im})_{1 \leq i \leq n ; 1\leq m\leq q} \in \mathcal{M}_{n \times q}$ and $V=(V_{i \ell})_{1 \leq i \leq n ; 1 \leq \ell \leq p} \in \mathcal{M}_{n \times p}$ respectively denote the vector of observations for individual $i$, the vector of all observations, the matrix of all individual parameters, and the matrix of all covariates. Let us also note $\theta=(\mu,\boldsymbol{\beta},\Gamma,\sigma^2)$ the unknown parameter, also-called the population parameter. 

The goal of the present work is to identify the relevant covariates, \textit{i.e.} those that best explain the variability between individuals. This can be framed as identifying the non-zero elements in $\boldsymbol{\beta}$. Indeed, for each $(\ell,m) \in \{1,\dots,p\} \times \{1,\dots,q\}$, coefficient $\beta_{\ell m}$ describes the influence of covariate $\ell$ on the individual values of the $m$-th parameter of the nonlinear curves $(\varphi_{im})_{1 \leq i \leq n}$. More precisely, $\beta_{\ell m}=0$ means that covariate $\ell$ has no effect on $(\varphi_{im})_{1 \leq i \leq n}$ whereas $\beta_{\ell m} \neq 0$ means that covariate $\ell$ gives some information on these parameters. Identifying the relevant covariates for all individual parameters amounts to selecting the support of $\boldsymbol{\beta}$, defined by $S^{*}$:
\begin{equation*}
    S^{*}=\bigg\{(\ell,m) \in \{1,\dots,p\} \times \{1,\dots,q\}  \bigg\vert \beta^{*}_{\ell m}\neq 0\bigg\},
\end{equation*}
where $\boldsymbol{\beta}^{*}$ is the true fixed effects matrix. To solve this problem in a high-dimensional context, that is when $p>>n$, it is natural to assume that each row of $\boldsymbol{\beta}^{*}$ is sparse, which means that many $\beta^{*}_{\ell m}$ are zero. 
An important point here is that model~\eqref{model} is a model with incomplete-data. 
Indeed, although the first layer~\eqref{niv1} is observed, it is not the case for the individual parameters $\boldsymbol{\varphi}$. The main difficulty here is that variable selection concerns latent variables of the model. 

\begin{remark}
    As in the application presented in Section \ref{Section_data_reel}, we could perform variable selection on $q'<q$ components of $\varphi_i$, with the choice of the $q'$ relevant components to be considered for variable selection driven by the initial biological question without changing the methodology presented in what follows.
\end{remark}

\subsection{Prior specification}
\label{Prior}

To solve this variable selection problem, it is convenient to adopt a Bayesian approach. The purpose of this section is to describe the prior distribution on $\theta=(\mu,\boldsymbol{\beta},\Gamma,\sigma^2)$. 
First, in order to find the non-zero coefficients of $\boldsymbol{\beta}$, a spike-and-slab mixture prior \citep{george, george97, rockova} is considered in a multivariate setting. To facilitate the formulation of this prior, binary latent variables $\boldsymbol{\delta}~=~(\delta_{\ell m})_{1\leq \ell \leq p ; 1 \leq m \leq q}$ are introduced, such as: 
\begin{equation}
\label{delta}
\forall 1 \leq \ell \leq p \text{ , } \forall 1 \leq m \leq q \text{ , }
    \delta_{\ell m}=\left\{
    \begin{array}{ll}
        1 & \mbox{if $(\ell,m)$ is to be included in model $S^{*}$ }, \\
        0 & \mbox{otherwise.}
    \end{array}
\right. 
\end{equation}
    
Thus, $\delta_{\ell m}=1$ indicates that the covariate $\ell$ provides information on the individual parameter $m$.  In other words, $\boldsymbol{\delta}$ characterises the support of $\boldsymbol{\beta}$. 
Note that contrary to many multivariate variable selection methods such as the one implemented in \textbf{glmnet} \citep{friedman}, the proposed procedure is not limited to selecting the same set of covariates for all dimensions $m$, i.e it is not required that $\forall 1 \leq m \leq q, \quad \delta_{\ell m} = \delta_{\ell}$. 
The support $S^{*}$ can therefore be reformulated as follows: 
\begin{equation}
\label{Sbeta}
    S^{*}=\bigg\{(\ell,m) \in \{1,\dots,p\} \times \{1,\dots,q\} \bigg \vert \delta^{*}_{\ell m}=1\bigg\},
\end{equation}
where $\boldsymbol{\delta}^{*}$ denotes the true support. Then, one would like to find $\widehat{\boldsymbol{\delta}}$ that maximises the posterior probability $\pi(\boldsymbol{\delta}\vert\boldsymbol{y})$, which corresponds to the most promising support. 

The prior formulations proposed here are based on the non-conjugate version of the hierarchical priors of \cite{george97}, summarised as follows: 
\begin{subequations}
\label{priors}
\begin{align}
    \pi(\beta_{\ell m} \vert \delta_{\ell m}) &= \mathcal{N}(0,(1-\delta_{\ell m})\nu_0 + \delta_{\ell m} \nu_1)\text{, } 1 \leq \ell \leq p \text{, } 1 \leq m \leq q\text{, } 0 < \nu_0 < \nu_1 \text{, } \label{prior_beta}\\
    \pi(\mu)&=  \mathcal{N}_q(0,\sigma^2_{\mu}I) \text{, with } \sigma^2_{\mu}>0 \text{, } \label{prior_mu}\\
    \pi(\sigma^2)&= \mathcal{IG}\left(\dfrac{\nu_{\sigma}}{2},\dfrac{\nu_{\sigma}\lambda_{\sigma}}{2}\right) \text{, with } \nu_{\sigma},\lambda_{\sigma}>0 \text{, } \label{prior_sigma} \\
    \pi(\Gamma)&= \mathcal{IW}\left(\Sigma_{\Gamma},d\right) \text{, with } \Sigma_{\Gamma}\in \mathcal{S}_{q}^{++}, d>0 \text{, } \label{prior_Gamma}\\
    \pi(\delta_{\ell m}\vert\alpha_{m})&= \alpha_{m}^{\delta_{\ell m}}(1-\alpha_{m})^{1-\delta_{\ell m}} \text{, with } \alpha_{m} \in [0,1], 1 \leq \ell \leq p, 1 \leq m \leq q \text{, } \label{prior_delta}\\
    \pi(\alpha_{m}) &= \textit{Beta}(a_m,b_m) \text{, with } a_m,b_m>0 \text{, } 1 \leq m \leq q.  \label{prior_alpha}
\end{align}
\end{subequations}

The key prior distribution used for variable selection in this method is the spike-and-slab Gaussian mixture prior~\eqref{prior_beta} on $\boldsymbol{\beta}$. In this prior, $\nu_0$ and $\nu_1$ are parameters controlling the penalisation inducing sparsity in the columns of $\boldsymbol{\beta}$. More precisely, $(\beta_{\ell m})_{1\leq \ell \leq p ; 1\leq m \leq q}$ are independent conditionally on $\boldsymbol{\delta}$, with $\pi(\beta_{\ell m}\vert\delta_{\ell m}~=~0)~=~\mathcal{N}(0,\nu_0)$ and $\pi(\beta_{\ell m}\vert\delta_{\ell m}~=~1)~=~\mathcal{N}(0,\nu_1)$. The general recommendation for this type of prior is to set $\nu_0$ small to encourage the exclusion of insignificant effects, and $\nu_1$ large enough to accommodate all plausible $\boldsymbol{\beta}$ values \citep[see][]{george97}. Indeed, when $\delta_{\ell m}=0$, the prior constrains $\beta_{\ell m}$ to very small values which implies that covariate $\ell$ has no impact in the individual parameter $m$ in the model. Thus, through the values $\delta_{\ell m}$, the spike-and-slab prior makes it possible to distinguish the selected covariates from the rest.

Note that, since $\boldsymbol{\varphi}$ is unobserved, one cannot simply centre the variable on which the selection is made as is usually the case in more standard models, and so the inclusion of an intercept $\mu$ is necessary. 
Thus, a vaguely informative Gaussian prior~\eqref{prior_mu} is used for $\mu$, with $\sigma^2_{\mu}$ large enough. 
This choice of prior has the advantage of simplifying the calculations for parameter inference, thanks to a useful reformulation $\boldsymbol{\Tilde{\beta}}~=~\begin{pmatrix}
\mu^\top \\
\boldsymbol{\beta} \\
\end{pmatrix}~\in~\mathcal{M}_{(p+1) \times q}$ and, for all $1~\leq~i~\leq~n$, $\Tilde{V}_{i}~=~(\Tilde{V}_{i\ell'})_{1 \leq \ell' \leq p+1}~=~(1,V_i)^\top~\in~\mathbb{R}^{p+1}$, so that $\mu + \boldsymbol{\beta}^\top V_i=\boldsymbol{\Tilde{\beta}}^\top \Tilde{V}_i$. Let $\Tilde{V}~=~(\Tilde{V}_{i\ell'})_{1 \leq i \leq n ; 1 \leq \ell' \leq p+1}$ such as $\boldsymbol{\varphi}=\Tilde{V} \boldsymbol{\Tilde{\beta}} + \xi$, where $\xi=\begin{pmatrix}
\xi_1^\top \\
\vdots \\
\xi_n^\top \\
\end{pmatrix} \in \mathcal{M}_{n \times q}$. Then, by introducing $\Tilde{\boldsymbol{\delta}}=(\Tilde{\delta}_{\ell' m})_{1 \leq \ell' \leq p+1 ; 1 \leq m \leq q }$ such as: \begin{equation*}
\forall 1 \leq \ell' \leq p+1 \text{ , } \forall 1 \leq m \leq q \text{ , }
    \Tilde{\delta}_{\ell' m}=\left\{
    \begin{array}{ll}
        1 & \mbox{if $\ell'=1$}, \\
        \delta_{\ell m} & \mbox{where $\ell=\ell'-1$, if $\ell' > 1$,}
    \end{array}
\right. 
\end{equation*}
to force the inclusion of the intercept in the model, Equations~\eqref{prior_beta} and \eqref{prior_mu} can be rewritten as: for $1 \leq \ell' \leq p+1$, $1 \leq m \leq q$,
\begin{equation}
\pi(\Tilde{\beta}_{\ell' m} \vert \Tilde{\delta}_{\ell' m}) = \mathcal{N}(0,(1-\Tilde{\delta}_{\ell' m})\nu_0 + \Tilde{\delta}_{\ell' m} (\mathds{1}_{\ell'>1}\nu_1 + \mathds{1}_{\ell'=1} \sigma_{\mu}^2)) 
\label{reformulation}
\end{equation}

For the variance parameter $\sigma^2$, an inverse-gamma prior is chosen \eqref{prior_sigma}, which prohibits negative values. One possibility is to set $\nu_{\sigma},\lambda_{\sigma}$ equal to $1$ for example, to make it relatively non-influential. For the inter-individual variance-covariance matrix $\Gamma$, the inverse-Wishart prior is chosen \eqref{prior_Gamma}, which is the multivariate extension of the inverse-Gamma density. The hyperparameter $\Sigma_{\Gamma}$ is a positive matrix which can be specified as $\Sigma_{\Gamma}=k I_q$ with $k$ chosen of comparable size to the maximum likely variance of $(\varphi_{im})_{1\leq i \leq m}$ for all $m\in \{1, \dots, q\}$. The hyperparameter $d$, which is a degree of freedom, can be chosen as the smallest integer value ensuring the existence of $\mathbb{E}[\Gamma]$, that is $q+2$.

Following \cite{rockova}, the i.i.d. Bernoulli prior~\eqref{prior_delta} is used for the inclusion variable $\boldsymbol{\delta}$, where the hyperparameters $(\alpha_m)_m$ can be seen as the proportion of relevant covariates for each individual parameter, and a Beta distribution prior~\eqref{prior_alpha} is chosen on each $\alpha_m$, for $m \in \{1, \dots, q\}$. To encourage sparsity in the model, \cite{castillo} suggest choosing $a_m$ small and $b_m$ large, for example $a_m=1$ and $b_m=p$ for all $m$. In the following, $\alpha$, $a$ and $b$  respectively denote the vectors $(\alpha_m)_{1 \leq m \leq q}$, $(a_m)_{1 \leq m \leq q}$ and $(b_m)_{1 \leq m \leq q}$. See \cite{rockova}, \cite{liquet2017bayesian}, \cite{deshpande2019simultaneous}, for more details on these choices of priors and the choice of hyperparameters values.

\section{Maximum \textit{a posteriori} inference and thresholding}
\label{inference}

The purpose of this section is to discuss the estimation of $\Theta=(\theta,\alpha)=(\boldsymbol{\Tilde{\beta}},\Gamma,\sigma^2,\alpha)$ in model~\eqref{model}~-~\eqref{priors}. Recall that $\Xi=(\nu_0, \nu_1, \sigma_{\mu}^2, \nu_{\sigma}, \lambda_{\sigma}, \Sigma_{\Gamma}, d, a, b)$ are fixed hyperparameters. 
In the following, model~\eqref{model}~-~\eqref{priors} is called SSNLME (Spike-and-Slab Non-Linear Mixed-Effects) model. Note that $\boldsymbol{\varphi}$, which is not observed, could be considered as a parameter to be estimated, included in $\Theta$. However, to design a scalable inference scheme, we consider it as a latent variable which we marginalise out of the posterior. This enables us to use an EM-type approach which is considerably faster than a full MCMC approach (see details in Subsection~\ref{comp_MCMC}). 

The EM-type approach proposed in Section~\ref{Model_select} requires to compute the maximum \textit{a posteriori} (MAP) estimator for $\Theta$: 
\begin{equation}
\label{MAP}
\widehat{\Theta}^{MAP}= \underset{\Theta \in \Lambda}{argmax} \quad \pi(\Theta\vert\boldsymbol{y})  \text{, with }  \pi(\Theta\vert\boldsymbol{y})=\dfrac{p_{\Theta}(\boldsymbol{y})\pi(\Theta)}{\int_{\Lambda} p_{\Theta}(\boldsymbol{y})\pi(\Theta) d\Theta} \text{ , }
\end{equation}
where $p_{\Theta}(\boldsymbol{y})$ and $\pi(\Theta)$ respectively denote the probability density of $\boldsymbol{y}$ conditionally to $\Theta$, and the prior density of $\Theta$, and $\Lambda$ denotes the parameter space. However, since the individual parameters $\boldsymbol{\varphi}$ are marginalised out, the $p_{\Theta}(\boldsymbol{y})$ distribution is not explicit. Denoting $Z=(\boldsymbol{\varphi},\boldsymbol{\delta})\in \mathcal{Z}$ the latent variables, the marginalised posterior distribution $\pi(\Theta\vert\boldsymbol{y})$ takes the form: 
$$\pi(\Theta\vert\boldsymbol{y})=\int_{\mathcal{Z}} \pi(\Theta,Z\vert\boldsymbol{y}) dZ \text{, with }  \pi(\Theta,Z\vert\boldsymbol{y})=\dfrac{p(\boldsymbol{y}\vert\Theta,Z)p(\Theta,Z)}{\int_{\mathcal{Z}} \int_{\Lambda} p(\boldsymbol{y}\vert\Theta,Z)p(\Theta,Z) d\Theta dZ } \text{ , } $$
where $p(\boldsymbol{y}\vert\Theta,Z)$ and $p(\Theta,Z)$ respectively designate the probability density of $\boldsymbol{y}$ conditionally to $(\Theta,Z)$, and the joint distribution of $\Theta$ and $Z$.
 
Targeting only the maximum \textit{a posteriori} replaces a sampling problem by an optimisation problem, which turns out to be much more scalable than exploring the full posterior.
Equation~\eqref{MAP} is an optimisation problem in an incomplete data model, which is gainfully tackled using the Stochastic Approximation version of the EM algorithm (SAEM, \cite{delyon}).
Often in the literature, the EM algorithm and its extensions are presented in the frequentist framework for the calculation of the maximum likelihood estimator (MLE). Nevertheless, these algorithms are also very well adapted to the computation of the MAP estimator \citep{dempster}.

\subsection{General description of SAEM algorithm}
\label{generality}

In this subsection, we consider the general framework of an incomplete data model with observations $\boldsymbol{y}$ and latent variables $Z$ that characterise the distribution of observations. 
It it assumed that the density of the complete data $(\boldsymbol{y},Z)$ is parameterised by $\Theta$, which is unknown and associated with a prior  $\pi(\Theta)$. 
The EM algorithm is iterative and allows to build a sequence $(\Theta^{(k)})_k$ of parameter estimates, which under certain regularity conditions converges to a local maximum of the observed posterior distribution $$\pi(\Theta\vert\boldsymbol{y})=\int \pi(\Theta,Z\vert\boldsymbol{y}) dZ, $$ (see \cite{delyon} for more details). However, this integral is generally intractable and the idea is to maximise it by iteratively maximising an easier quantity: $$Q(\Theta\vert\Theta')~=~\mathbb{E}_{Z\vert(\boldsymbol{y},\Theta')}[\log(\pi(\Theta,Z\vert\boldsymbol{y}))\vert\boldsymbol{y},\Theta'], $$ 
the conditional expectation of the complete log-posterior $\log(\pi(\Theta,Z\vert\boldsymbol{y}))$ given the observations $\boldsymbol{y}$ and the current value of the parameter estimates $\Theta'$. However, the quantity $Q(\Theta\vert\Theta')$ does not always have a closed form. This is especially the case in non-linear mixed-effects models like SSNLME model. However, even though this expectation cannot be computed in closed-form, it can be approximated by simulation. One solution proposed by \cite{wei} is to replace the E-step, \textit{i.e.} the computation of the $Q$ quantity, by a Monte-Carlo approximation based on a large number of independent simulations of the latent variables $Z$, called the MCEM algorithm. However, this large number of simulations is computationally expensive. The SAEM algorithm is an other alternative that replaces the E-step by a stochastic approximation based on a single simulation of the latent variables, which considerably reduces the computational cost \citep{delyon}. More precisely, the E-step of the EM algorithm is replaced by two steps: a simulation step (S-step) and a stochastic approximation step (SA-step). Then the $k$-th iteration of the SAEM algorithm proceeds as follows: 
\begin{enumerate}
    \item \textbf{S-step (Simulation):} simulate a realisation $Z^{(k)}$ of the latent variables according to the conditional distribution $\pi(Z\vert\boldsymbol{y},\Theta^{(k)})$.
    \item \textbf{SA-step (Stochastic Approximation):} update the approximation $Q_{k+1}(\Theta)$ of $Q(\Theta\vert\Theta^{(k)})$ by a stochastic approximation method, according to: $$Q_{k+1}(\Theta)=Q_k(\Theta) + \gamma_k(\log \pi(\Theta,Z^{(k)}\vert\boldsymbol{y}) - Q_k(\Theta)),$$
    where $(\gamma_k)_k$ is a sequence of step sizes decreasing towards 0 such that $\forall k$, $\gamma_k\in [0,1]$, $\sum_k \gamma_k=\infty$ and $\sum_k \gamma_k^2<\infty$.
    \item \textbf{M-step (Maximisation):} update the parameter value by computing: $$\Theta^{(k+1)}=\underset{\Theta \in \Lambda}{\operatorname{argmax}} \text{ } Q_{k+1}(\Theta).$$
\end{enumerate}

\begin{remark}
\label{expon_family}
If the model belongs to the curved exponential family, that is the complete log-posterior can be written as: $$\log(\pi(\Theta,Z\vert\boldsymbol{y}))=-\psi(\Theta)+ \bigg\langle S(\boldsymbol{y},Z), \phi(\Theta) \bigg\rangle, $$
where $\psi$ and $\phi$ denote two functions of $\Theta$, with $\langle \cdot , \cdot \rangle$ denoting the scalar product, and $S(\boldsymbol{y},Z)$ the minimal sufficient statistics of the model, then, $$Q(\Theta\vert\Theta^{(k)})=-\psi(\Theta)+\bigg\langle \mathbb{E}_{Z\vert(\boldsymbol{y},\Theta^{(k)})}[S(\boldsymbol{y},Z)\vert\boldsymbol{y},\Theta^{(k)}], \phi(\Theta) \bigg\rangle.$$ 
It is therefore sufficient to focus on the minimal sufficient statistics instead of $Q(\Theta\vert\Theta^{(k)})$ itself. More precisely, the SA-step and M-step of the SAEM algorithm are replaced by: 
\begin{itemize}
    \item \textbf{SA-step:} update $S_{k+1}$, the stochastic approximation of $\mathbb{E}_{Z\vert(\boldsymbol{y},\Theta^{(k)})}[S(\boldsymbol{y},Z)\vert\boldsymbol{y},\Theta^{(k)}]$, according to: $$S_{k+1}=S_k + \gamma_k(S(\boldsymbol{y},Z^{(k)})-S_k).$$
    \item \textbf{M-step:} update the parameter value by computing: $$\Theta^{(k+1)}=\underset{\Theta \in \Lambda}{\operatorname{argmax}} \text{ } \Big\{ -\psi(\Theta) + \langle S_{k+1}, \phi(\Theta) \rangle \Big\}.$$
\end{itemize}
Note that theoretical convergence results of the SAEM algorithm are provided in \cite{delyon} under the assumption that the model belongs to the curved exponential family. The SSNLME model belongs to this curved exponential family. 
\end{remark}

Note that the simulation step is not always directly feasible. This is particularly true in non-linear mixed-effects models since the conditional distribution of the latent variables knowing the observations and the current value of the parameters is known only to a nearest multiplicative constant. \cite{kuhn} proposed an alternative which consists in coupling the SAEM method with an MCMC procedure. Interestingly, convergence of the MCMC part at each S-step is not necessary and only a few MCMC iterations are required in practice (see \cite{kuhn} and \cite{kuhn2005} for practical and theoretical considerations about this algorithm). In the following, this extension is called MCMC-SAEM.

\subsection{Central decomposition of the Q quantity in spike-and-slab non-linear mixed-effects models}
\label{Decomposition}

In the following, the notations from Section~\ref{stat_model} are used again, and $\vert\vert \cdot \vert\vert$ denotes the Euclidean norm on $\mathbb{R}^n$. The SSNLME model~\eqref{model}~-~\eqref{priors} is a particular latent variables model with $\boldsymbol{y}=(y_{ij})_{i,j}$ and $Z=(\boldsymbol{\varphi},\boldsymbol{\delta})$. The aim here is to decompose the $Q$ quantity of the SAEM algorithm in the particular case of the SSNLME model, allowing to describe an algorithm for computing the MAP estimator of $\Theta$ in the following subsection.

First, note that, by using the tower property of conditional expectation, the quantity $Q(\Theta\vert\Theta^{(k)})$ in model~\eqref{model}~-~\eqref{priors} is written as:
\begin{align*}
    Q(\Theta\vert\Theta^{(k)}) &= \mathbb{E}_{(\boldsymbol{\varphi},\boldsymbol{\delta}) \vert(\boldsymbol{y},\Theta^{(k)})}[\log(\pi(\Theta,\boldsymbol{\varphi},\boldsymbol{\delta}\vert\boldsymbol{y}))\vert\boldsymbol{y},\Theta^{(k)}] \\
    &= \mathbb{E}_{\boldsymbol{\varphi}\vert(\boldsymbol{y},\Theta^{(k)})}\left[\overset{\sim}{Q}(\boldsymbol{y},\boldsymbol{\varphi},\Theta,\Theta^{(k)})\bigg\vert\boldsymbol{y},\Theta^{(k)}\right],
\end{align*}
    
where 
\begin{equation}
    \label{Qtilde}
    \overset{\sim}{Q}(\boldsymbol{y},\boldsymbol{\varphi},\Theta,\Theta^{(k)})= \mathbb{E}_{\boldsymbol{\delta}\vert(\boldsymbol{\varphi},\boldsymbol{y},\Theta^{(k)})}[\log(\pi(\Theta,\boldsymbol{\varphi},\boldsymbol{\delta}\vert\boldsymbol{y})) \vert \boldsymbol{\varphi},\boldsymbol{y},\Theta^{(k)}].
\end{equation}
It is interesting to write $Q(\Theta\vert\Theta^{(k)})$ like this because $\overset{\sim}{Q}(\boldsymbol{y},\boldsymbol{\varphi},\Theta,\Theta^{(k)})$ has a closed form. Indeed, Proposition~\ref{PropQtilde} shows that $\overset{\sim}{Q}(\boldsymbol{y},\boldsymbol{\varphi},\Theta,\Theta^{(k)})$ can be decomposed such as: \begin{equation}
\label{separability}
    \overset{\sim}{Q}(\boldsymbol{y},\boldsymbol{\varphi},\Theta,\Theta^{(k)})= C  + \overset{\sim}{Q}_{1}(\boldsymbol{y},\boldsymbol{\varphi},\theta, \Theta^{(k)}) + \overset{\sim}{Q}_{2}(\alpha,\Theta^{(k)}),
\end{equation}
where $\overset{\sim}{Q}_{1}$ and $\overset{\sim}{Q}_{2}$ have a closed form given in Proposition~\ref{PropQtilde}. Thus, the separability of \eqref{separability} into two distinct functions, $\overset{\sim}{Q}_{1}$ which depends on $(\boldsymbol{y},\boldsymbol{\varphi},\theta, \Theta^{(k)})$ and $\overset{\sim}{Q}_{2}$ on $(\alpha,\Theta^{(k)})$, allows to update the estimations of $\theta$ and $\alpha$ independently from one another. Moreover, since $\overset{\sim}{Q}_{2}$ does not depend on $\boldsymbol{\varphi}$, Proposition~\ref{PropQtilde} allows to write that: 
\begin{equation}
\label{Q_simple}
    Q(\Theta\vert\Theta^{(k)})=C + \mathbb{E}_{\boldsymbol{\varphi}\vert(\boldsymbol{y},\Theta^{(k)})}\left[\overset{\sim}{Q}_{1}(\boldsymbol{y},\boldsymbol{\varphi},\theta,\Theta^{(k)})\bigg\vert\boldsymbol{y},\Theta^{(k)}\right] + \overset{\sim}{Q}_{2}(\alpha,\Theta^{(k)}).
\end{equation}
However, even if $\overset{\sim}{Q}(\boldsymbol{y},\boldsymbol{\varphi},\Theta,\Theta^{(k)})$ has a closed form, this is not the case of $Q(\Theta\vert\Theta^{(k)})$ because the function $g$ is non-linear with respect to $\varphi_i$, and so $\pi(\boldsymbol{\varphi}\vert\boldsymbol{y},\Theta^{(k)})$ is only known to a nearest multiplicative constant. Thus, it is necessary to use a stochastic approximation method to approximate $\mathbb{E}_{\boldsymbol{\varphi}\vert(\boldsymbol{y},\Theta^{(k)})}\left[\overset{\sim}{Q}_{1}(\boldsymbol{y},\boldsymbol{\varphi},\theta,\Theta^{(k)})\bigg\vert\boldsymbol{y},\Theta^{(k)}\right]$ in Equation~\eqref{Q_simple}. The originality of the present extension of the MCMC-SAEM algorithm is that it combines an exact computation $\overset{\sim}{Q}_{2}(\alpha,\Theta^{(k)})$ and a stochastic approximation of $\mathbb{E}_{\boldsymbol{\varphi}\vert(\boldsymbol{y},\Theta^{(k)})}\left[\overset{\sim}{Q}_{1}(\boldsymbol{y},\boldsymbol{\varphi},\theta,\Theta^{(k)})\bigg\vert\boldsymbol{y},\Theta^{(k)}\right]$ instead of a stochastic approximation of the entire quantity $Q(\Theta\vert\Theta^{(k)})$. This results in the combination of an exact EM algorithm and of an MCMC-SAEM algorithm for the estimation of $\alpha$ and $\theta$ respectively.

Also, let us notice that $\overset{\sim}{Q}_{1}(\boldsymbol{y},\boldsymbol{\varphi},\theta,\Theta^{(k)})$ takes an exponential form. 
Thus, according to Remark~\ref{expon_family}, it suffices to approximate stochastically $\mathbb{E}_{\boldsymbol{\varphi}\vert(\boldsymbol{y},\Theta^{(k)})}\left[S(\boldsymbol{y},\boldsymbol{\varphi})\big\vert\boldsymbol{y},\Theta^{(k)}\right]$ at SA-step (see Appendix \ref{app_MCMCSAEM} for more details about $S(\boldsymbol{y},\boldsymbol{\varphi})$ and this exponential form). 

Algorithm~\ref{algoMCMCSAEM} in Appendix \ref{app_MCMCSAEM} summarises the proposed extension of the MCMC-SAEM algorithm for computing the MAP estimator of $\Theta$ in the SSNLME model.

\begin{proposition}
\label{PropQtilde}
Consider $\overset{\sim}{Q}(\boldsymbol{y},\boldsymbol{\varphi},\Theta,\Theta^{(k)})$ defined by Equation~\eqref{Qtilde} where \\ $\Theta=(\boldsymbol{\Tilde{\beta}},\Gamma,\sigma^2,\alpha)$. Then:
\begin{equation*}
    \overset{\sim}{Q}(\boldsymbol{y},\boldsymbol{\varphi},\Theta,\Theta^{(k)})= C  + \overset{\sim}{Q}_{1}(\boldsymbol{y},\boldsymbol{\varphi},\theta, \Theta^{(k)}) + \overset{\sim}{Q}_{2}(\alpha,\Theta^{(k)}),
\end{equation*}
where $C$ is a normalisation constant which does not depend on $\Theta$, and with: 
\begin{align*}
    \overset{\sim}{Q}_{1}(\boldsymbol{y},\boldsymbol{\varphi},\theta, \Theta^{(k)})= &-\dfrac{1}{2\sigma^2}\sum_{i=1}^n \sum_{j=1}^{n_i} (y_{ij}-g(\varphi_i,t_{ij}))^2 - \dfrac{n_{tot} + \nu_{\sigma}+2}{2}\log(\sigma^2) - \dfrac{\nu_{\sigma} \lambda_{\sigma}}{2\sigma^2} \\
    &- \dfrac{1}{2} \mathrm{Tr} \left((\boldsymbol{\varphi} - \Tilde{V}\boldsymbol{\Tilde{\beta}})^\top (\boldsymbol{\varphi} - \Tilde{V}\boldsymbol{\Tilde{\beta}}) \Gamma^{-1}\right)  - \dfrac{1}{2}\sum_{m=1}^q\sum_{\ell'=1}^{p+1} \Tilde{\beta}_{\ell' m}^2 \overset{\sim}{d}_{\ell' m}^{*}(\Theta^{(k)}) \\
    &- \dfrac{n+d+q+1}{2}\log(\vert\Gamma\vert) -\dfrac{1}{2} \mathrm{Tr}(\Sigma_{\Gamma} \Gamma^{-1})
\end{align*}
where $\vert A\vert$ and $\mathrm{Tr}(A)$ respectively denote the determinant and the trace of a matrix $A$, and 
\begin{align*}
    \overset{\sim}{Q}_{2}(\alpha,\Theta^{(k)})=&\sum_{m=1}^q \log\left(\sqrt{\dfrac{\nu_0}{\nu_1}} \dfrac{\alpha_m}{1-\alpha_m}\right) \sum_{\ell=1}^{p}  p_{\ell m}^{*}(\Theta^{(k)})  + \\
    &(a_m-1)\log(\alpha_m) + (p+b_m-1)\log(1-\alpha_m).
\end{align*}
Quantities $p_{\ell m}^{*}(\Theta^{(k)})$, $1\leq \ell \leq p$, $1 \leq m \leq q$, and $\overset{\sim}{d}_{\ell' m}^{*}(\Theta^{(k)})$, $1 \leq \ell' \leq p+1$, $1 \leq m \leq q$, are defined as follows:
\begin{equation}
\label{p_l}
    p_{\ell m}^{*}(\Theta^{(k)})=\mathbb{E}[\delta_{\ell m}\vert\boldsymbol{\varphi},\boldsymbol{y},\Theta^{(k)}]= \dfrac{\alpha_m^{(k)}\phi_{\nu_1}(\beta_{\ell m}^{(k)})}{\alpha_m^{(k)}\phi_{\nu_1}(\beta_{\ell m}^{(k)})+ (1-\alpha_m^{(k)})\phi_{\nu_0}(\beta_{\ell m}^{(k)})}
\end{equation}
where $\phi_{\nu}(\cdot)$ is the normal density with zero mean and variance $\nu$, and 
\begin{align}
            \overset{\sim}{d}_{\ell' m}^{*}(\Theta^{(k)})&=\mathbb{E}\left[\dfrac{1}{(1-\Tilde{\delta}_{\ell' m})\nu_0 + \Tilde{\delta}_{\ell' m} (\mathds{1}_{\ell'>1}\nu_1 + \mathds{1}_{\ell'=1} \sigma_{\mu}^2)}\bigg\vert\boldsymbol{\varphi},\boldsymbol{y},\Theta^{(k)}\right] \notag \\
    &=\dfrac{1}{\sigma^2_{\mu}}\mathds{1}_{\ell'=1}+ \left(\dfrac{1-p_{\ell m}^{*}(\Theta^{(k)})}{\nu_0}+\dfrac{p_{\ell m}^{*}(\Theta^{(k)})}{\nu_1}\right)\mathds{1}_{\ell'\neq 1} \label{d_l}
\end{align}
where $\ell = \ell'-1$.
\end{proposition}

\begin{remark}
Note that $\mathbb{E}[\delta_{\ell m}\vert\boldsymbol{\varphi},\boldsymbol{y},\Theta^{(k)}]=\mathbb{E}[\delta_{\ell m}\vert\Theta^{(k)}]$ because the posterior distribution of $\boldsymbol{\delta}$ given $(\boldsymbol{\varphi},\boldsymbol{y},\Theta^{(k)})$ depends on $\boldsymbol{y}$ and $\boldsymbol{\varphi}$ only through the current estimates $\Theta^{(k)}$. 
\end{remark}

\begin{remark}
Note that for a linear mixed-effects model, that is when $g$ is linear with respect to $\varphi_i$, a classical EM algorithm is applicable.
\end{remark} 
\color{black}

\subsection{Estimator thresholding}
\label{Thresholding}

As in \cite{rockova}, after obtaining an estimator $\widehat{\Theta}^{MAP}$, the support $S^{*}$, defined in Equation~\eqref{Sbeta}, can be naturally estimated as the most probable model conditionally on $\widehat{\Theta}^{MAP}$. Indeed, for all $m~\in~\{1, \dots, q\}$ and for all $\ell \in \{1,\dots,p\}$, since $\pi(\delta_{\ell m}=1\vert\widehat{\alpha}_m^{MAP})=\widehat{\alpha}_m^{MAP}$ and $\pi(\delta_{\ell m}=0\vert\widehat{\alpha}_m^{MAP})=1-\widehat{\alpha}_m^{MAP}$, the \textit{a posteriori} inclusion probability of the covariate $\ell$ knowing $\widehat{\Theta}^{MAP}$ for the individual parameter $m$ can be obtained as: 
\begin{equation*}
    \mathbb{P}(\delta_{\ell m}=1 \vert\boldsymbol{y}, \widehat{\beta}^{MAP}_{\ell m}, \widehat{\alpha}_m^{MAP})=\dfrac{\pi_1(\widehat{\beta}^{MAP}_{\ell m})\widehat{\alpha}_m^{MAP}}{\pi_1(\widehat{\beta}^{MAP}_{\ell m})\widehat{\alpha}_m^{MAP} + \pi_0(\widehat{\beta}^{MAP}_{\ell m})(1-\widehat{\alpha}_m^{MAP})},
\end{equation*}
where $\pi_k(\widehat{\beta}^{MAP}_{\ell m})=\pi(\widehat{\beta}^{MAP}_{\ell m}\vert\delta_{\ell m}=k)$ for $k\in\{0,1\}$.
Then, $\widehat{\boldsymbol{\delta}}$, which is the most probable $\boldsymbol{\delta}$ knowing that $\Theta=\widehat{\Theta}^{MAP}$, can be computed as follows: 
\begin{align*}
    \widehat{\delta}_{\ell m}=1 &\Longleftrightarrow \mathbb{P}(\delta_{\ell m}=1 \vert\boldsymbol{y}, \widehat{\beta}^{MAP}_{\ell m},\widehat{\alpha}_m^{MAP}) \geq  0.5 \\
    &\Longleftrightarrow \vert\widehat{\beta}^{MAP}_{\ell m}\vert \geq \sqrt{2\dfrac{\nu_0 \nu_1}{\nu_1-\nu_0} \log\left( \sqrt{\dfrac{\nu_1}{\nu_0}} \dfrac{1-\widehat{\alpha}_m^{MAP}}{\widehat{\alpha}_m^{MAP}}\right)} = s_{\beta}(\nu_0,\nu_1,\widehat{\alpha}_m^{MAP}).
\end{align*}
Note that this estimator can be seen as a local version of the median probability model of \cite{barbieri2004optimal}.
Thus, the following subset of covariates for the individual parameter $m$ is selected via a thresholding operation: 
\begin{equation}
\label{S_m}
    \widehat{S}=\bigg\{(\ell,m) \in \{1,\dots,p\} \times \{1,\dots,q\} \text{ } \bigg\vert \text{ } \vert\widehat{\beta}^{MAP}_{\ell m}\vert \geq s_{\beta}(\nu_0,\nu_1,\widehat{\alpha}_m^{MAP}) \bigg\}.
\end{equation}

\begin{remark}
Note that threshold $s_{\beta}(\nu_0,\nu_1,\widehat{\alpha}^{MAP}_{m})$ is the same for all the covariates but depends on the individual parameter $m$ and on the values of the spike and slab hyperparameters $\nu_0$ and $\nu_1$ which act as tuning parameters for the penalty.
\end{remark}

\begin{remark}
It is interesting to note that the thresholding rule is unchanged from the easier situation where the individual parameters $\varphi_i$'s would have been directly observed, which would have fit into the framework treated in \cite{rockova}.
\end{remark}

\section{Covariate selection procedure}
\label{Model_select}

Similarly to Lasso regression \citep{tibshirani_regression_1996}, it is interesting to exploit the flexibility of the spike-and-slab prior to study different levels of sparsity in the $\boldsymbol{\beta}$'s columns, and thanks to the speed of the MCMC-SAEM algorithm, it is possible to explore a grid of values for the spike hyperparameter $\nu_0$ rather than focusing on a single value. Indeed, mechanically, the higher $\nu_0$ is, the fewer covariates are included in the estimated support of $\boldsymbol{\beta}$'s columns. This is why it is more interesting to look at a grid of values and then use a model selection criterion to choose the optimal model. Let us denote $\Delta$ this grid, and $\vert\Delta\vert$ the number of grid points. Then, for all $\nu_0 \in \Delta$, the MCMC-SAEM algorithm is executed to obtain the MAP estimate of $\Theta$, $\widehat{\Theta}^{MAP}_{\nu_0}$, which is then used to determine a subset of relevant covariates for all individual parameters, $\widehat{S}_{\nu_0}$, as explained by Equation~\eqref{S_m} in Subsection~\ref{Thresholding}. This first step reduces the total collection of $2^{pq}$ possible models to a smaller collection of $\vert\Delta\vert \ll  2^{pq}$ promising sub-models $(\widehat{S}_{\nu_0})_{\nu_0 \in \Delta}$ with high posterior probability. Next, a model selection criterion can be applied to choose the "best" model from this collection.

As explained in \cite{rockova}, a possible criterion is to maximise, along the grid, the marginal posterior of $\boldsymbol{\delta}$ under the prior with $\nu_0=0$. This corresponds to the so-called Dirac-and-slab prior, where the spike is a Dirac distribution \citep{mitchell}. However, in our case, it is not possible to have an explicit expression for this marginal and it is also difficult to obtain it numerically, so this criterion is not convenient. 

However, as the collection of models has been reduced to a small sub-collection $(\widehat{S}_{\nu_0})_{\nu_0 \in \Delta}$, that contains at most $\vert\Delta\vert$ models, an information criterion can be used to choose the final model. Thus, covariate selection would consist in choosing the "best" $\nu_0 \in \Delta$, that is noted $\hat \nu_0$, as: 
\begin{equation}
\label{nu0_chap}
    \hat\nu_0 = \underset{\nu_0 \in \Delta}{\text{argmin}} \Big \{ \text{crit}(\widehat{S}_{\nu_0} ) \Big\},
\end{equation}
where
\begin{equation}
\label{crit}
   \text{crit}(\widehat{S}_{\nu_0}) = -2 \log \left(p(\boldsymbol{y};\hat\theta_{\nu_0})\right) + \text{pen}(\nu_0),
\end{equation}
with:
\begin{itemize}
    \item $\log \left(p(\boldsymbol{y};\theta)\right)$ the log-likelihood of model~\eqref{model},
    \item $\hat\theta_{\nu_0}=(\widehat{\boldsymbol{\Tilde{\beta}}}_{\nu_0},\widehat{\Gamma}_{\nu_0},\widehat{\sigma}_{\nu_0}^{2})$ a point estimator of the parameter $\theta~=~(\boldsymbol{\Tilde{\beta}},\Gamma,\sigma^2)$ in sub-model $\widehat{S}_{\nu_0}$,
    \item $\text{pen}(\nu_0)$ a penalty function which penalizes the complexity of $\widehat{S}_{\nu_0}$.
\end{itemize}

There are many ways to define the penalty in the criterion~\eqref{crit}, e.g. AIC \citep{akaike1998information}, BIC \citep{schwarz, delattre}, DIC \citep{spiegelhalter2002bayesian}, etc. 
Here, a pragmatic and effective choice is to use the eBIC (extended Bayesian Information Criterion, \cite{chen}) which is tailored to the high-dimensional setting. 
Indeed, eBIC's penalty has the following form: 

\begin{equation}
\label{eBIC}
    \text{pen}_{\text{eBIC}}(\nu_0)=\lvert \widehat{S}_{\nu_0} \rvert\times \log(n)+ 2\log\left(\binom{pq}{\lvert \widehat{S}_{\nu_0} \rvert}\right),
\end{equation}

where $\lvert \widehat{S}_{\nu_0} \rvert$ is the size of this support, which allows to take into account that the number of possible models with $r \leq pq$ covariates increases quickly as $r$ increases. The eBIC uses $\hat\theta_{\nu_0}$ the maximum likelihood estimate (MLE).
Note that this MLE and log-likelihood that are required to compute criteria with a penalty of the form~\eqref{eBIC} do not have an explicit form here because the individual parameters are latent and the function $g$ is non-linear with respect to $\varphi_i$. They are calculated using an MCMC-SAEM algorithm and importance sampling techniques respectively (see \textit{e.g.} \cite{kuhn2005} and \cite{lavielle} for details). Note that some Bayesian model selection criteria could also be appropriate instead of eBIC.

The proposed variable selection procedure can be summarised as in Algorithm~\ref{algoSAEMVS} in Appendix \ref{app_SAEMVS}. In the following, this procedure is called SAEMVS (SAEM Variable Selection). A detailed example of the application of this procedure on a toy example is presented in a supporting web material given in the section "Data and code availability".

\begin{remark}
Note that for Algorithm~\ref{algoSAEMVS} it is possible to parallelise the computations along the grid because the outputs of the algorithm for two given values of $\nu_0 \in \Delta$ do not depend on each other. 
\end{remark}

\section{Numerical experiments}
\label{simul}

This section studies the performance of the global selection procedure, named SAEMVS, in the non-linear mixed effects model with spike-and-slab prior, called the SSNLME model. As a reminder, this procedure is summarised in Algorithm~\ref{algoSAEMVS} in Appendix \ref{app_SAEMVS}, and uses an MCMC-SAEM algorithm for inference, which is detailed in Algorithm~\ref{algoMCMCSAEM} in Appendix \ref{app_MCMCSAEM}.

The numerical study is divided into three parts. The first part is a comparison with strategies that can be easily implemented from existing methods.  It aims essentially at demonstrating the interest of carrying out the selection of covariates from the data of all the individuals simultaneously thanks to the mixed effects model. The second part studies in great detail the influences of the number of subjects, the number of covariates, the signal-to-variability ratio and the collinearity between covariates on the performance of SAEMVS. To both show the flexibility of our approach and simplify the presentation of this study, this second part is conducted on another nonlinear mixed effects model, this time with a one-dimensional random effect. In the third part, a comparison of SAEMVS in terms of computation time with an MCMC implementation is presented, to quantify precisely the speed improvement afforded by the SAEM algorithm.

\subsection{Comparison with a two-step approach}
\label{simul_multi}

The proposed approach is compared to a standard solution readily available to a practitioner using existing tools. The goal being to study the impact of covariates on specific parameters of the nonlinear models, a manageable approach would be to proceed in two steps and fit independently the nonlinear model to each individual, then perform variable selection using as dependent variables the estimated parameters. This second step can be carried out using for instance the popular \textbf{glmnet} package, which allows multivariate response variable selection. This strategy is expected to work fine in data-rich scenarios when each parameter can be estimated very precisely, but it loses the uncertainty on the estimated parameters and the shrinkage property of the mixed-effect model. We are not aware of alternative solutions with freely available code that a practitioner could use for a high-dimensional problem as discussed in the introduction. The interest of SAEMVS, which considers the mixed effect model and fully embraces the uncertainty in the estimations, is shown in the context of repeated measures. 
We show that SAEMVS performs better than the two-steps approach on two grounds:
\begin{itemize}
    \item when the estimation of the individual parameters is more challenging, for instance in the classical scenario of lost to follow-up patients where some observation times may be missing,
    \item when different sets of covariates have an impact on the different dimensions of the response.
\end{itemize}

\subsubsection{Model and simulation design}
The following model, commonly used in pharmacokinetics, is considered:
\begin{equation}
\label{model_comp}
    \left\{
    \begin{array}{ll}
        y_{ij} = \dfrac{D \varphi_{i1}}{\mathcal{V}} \varphi_{i1} - \varphi_{i2} \left(\mathrm{e}^{-\dfrac{\varphi_{i2}}{\mathcal{V}}} t_{ij} - \mathrm{e}^{- \varphi_{i1} t_{ij}} \right) + \varepsilon_{ij} ,& \varepsilon_{ij} \overset{\text{i.i.d.}}{\sim} \mathcal{N}(0,\sigma^2), \\
        \varphi_i = \mu + \boldsymbol{\beta}^\top V_i+ \xi_i, & \xi_i \overset{\text{i.i.d.}}{\sim} \mathcal{N}_q(0,\Gamma),
    \end{array}
\right.
\end{equation}
where $\varphi_i=(\varphi_{i1},\varphi_{i2})^\top$, namely $q=2$. The constants $D$ and $\mathcal{V}$ are set to $100$ and $30$ respectively. The aim is to obtain a set of active covariates for the two individual parameters. For this study, the parameters are set to: $n=200$ individuals, $p=500$ covariates, $\sigma^2=10^{-3}$, $\Gamma= \left (
   \begin{array}{cc}
      0.2 & 0.05 \\
      0.05 & 0.1 \\ 
   \end{array}
\right)$ so that the correlation between the two individual parameters is of the order of $0.35$, $\mu=(6,8)^\top$, and $\boldsymbol{\beta}= \left (
   \begin{array}{cccccccc}
      3 & 2 & 1 & 0 & 0 & 0 & \dots & 0 \\
      0 & 0 & 3 & 2 & 1 & 0 & \dots & 0 
   \end{array}
\right)^\top$. The individual covariates $V_i \in \mathbb{R}^p$, $1 \leq i \leq n$, are simulated independently according to a binomial distribution with a success probability of $0.2$. Note that the support of the $\boldsymbol{\beta}$ is not the same for each dimension. It seems unrealistic that the covariates impacting $\varphi_{i1}$ would be exactly the same as those impacting $\varphi_{i2}$. The covariates are standardised. Thus, $100$ data-sets are simulated according to these parameter values where the number of observations by individual is fixed to $n_1 = \dots = n_n = 12$ with the following observation time points: $(t_{i1}, \dots , t_{i12})=(0.05,0.15,0.25,0.4,0.5,0.8,1,2,7,12,24,40)$.

This comparison was carried out on different scenarios, one corresponding to a data-rich situation where the two-step approach is expected to work relatively well, and a more challenging scenario where it is expected that there will be a difference between the two-step method and SAEMVS. The two scenarios correspond to different observation periods for each individual:

\begin{enumerate}
    \item \textbf{Complete data-set.} This is the baseline scenario where all individuals are observed during the entire experiment. 
    \item \textbf{Partial observations.} For each $p_{\textrm{partial}} \in \{0.1, 0.2, 0.3, 0.4 \}$, the other scenarios correspond to the case where $N_1=p_{\textrm{partial}} n$ individuals are assumed to be no longer part of the experiment after the 3rd observation time, that is: for each previously simulated data-set, only the first 3 observation times are kept for the first $N_1$ individuals, and all observation times for the remaining $N_2=n-N_1$ individuals.
\end{enumerate}

\subsubsection{Competing methods}

For this comparison study, SAEMVS is compared to two other procedures. For both of these procedures, the first step consists in estimating the $\varphi_i$'s individual-by-individual with the method of least squares thanks to the \textbf{nlm} R function (Non-Linear Minimization, see \cite{schnabel1985modular}). Then, the Lasso method from the \textbf{glmnet} R package \citep{friedman} is applied in its multivariate and univariate versions on the second level of the model \eqref{model_comp} using the estimated $\varphi_i$'s:
\begin{itemize}
    \item \textbf{Multivariate setting.} To take into account the correlation between the two individual parameters, the multi-response Gaussian family is used for the \textbf{glmnet} function on the estimated $\varphi_i$'s. Note that this function uses a group Lasso penalty, forcing the two individual parameters to have the same support. This method is called "mgaussian" in the following.
    \item \textbf{Univariate setting.} To circumvent this constraint, we also compare SAEMVS to the case where selection is made on the two individual parameters independently. For this, the \textbf{glmnet} function with the Gaussian family is used for each of the individual parameters separately on the estimated $\varphi_{i1}$'s and $\varphi_{i2}$'s. This method is called "gaussian" in the following.
\end{itemize}
For both of these methods, a cross-validation procedure is performed using the \textbf{cv.glmnet} function and the largest $\lambda$ at which the mean squared error (MSE) is within one standard error of the smallest MSE is chosen as recommended in \cite{hastie2009elements}. The algorithmic settings of SAEMVS are given in Appendix \ref{setting_SAEMVS_multi} with an example of convergence graphs of the MCMC-SAEM algorithm.

\subsubsection{Results }

First, Table \ref{Tab_phi_chap} compares the mean estimation errors for the two individual parameters using the first step of the two-step approach. It should be noted that 3 observation times appear sufficient to estimate the first parameter accurately, but insufficient for the second, as can be seen with the increasing estimation error for decreasing amount of data in Table \ref{Tab_phi_chap}.

\begin{table*}[!ht]
\begin{minipage}{350pt}
\caption{Comparison of the mean estimation errors for the first individual parameter (MEE1) and the second (MEE2) calculated on all individuals over the 100 data-sets using the first step of the two-step approach.}\label{Tab_phi_chap}%
\centering
\begin{tabular}{@{}llllll@{}}
\toprule
$p_{\textrm{partial}}$ & 0 & 0.1 & 0.2 & 0.3 & 0.4\\
\midrule
\textbf{MEE1}\footnote{MEE1 is the mean of the difference in absolute value between the true $\varphi_{i1}$ and its estimate over all the individuals and the 100 data-sets.}  & 0.088 & 0.10 & 0.10 & 0.11 & 0.12  \\
\textbf{MEE2}\footnote{MEE2 is the mean of the difference in absolute value between the true $\varphi_{i2}$ and its estimate over all the individuals and the 100 data-sets.}  & 0.12 & 0.62 & 1.14 & 1.68 & 2.20  \\
\toprule
\end{tabular}
\end{minipage}
\end{table*}

\begin{figure}[!ht]
    \centering
    \includegraphics[width=0.9\textwidth]{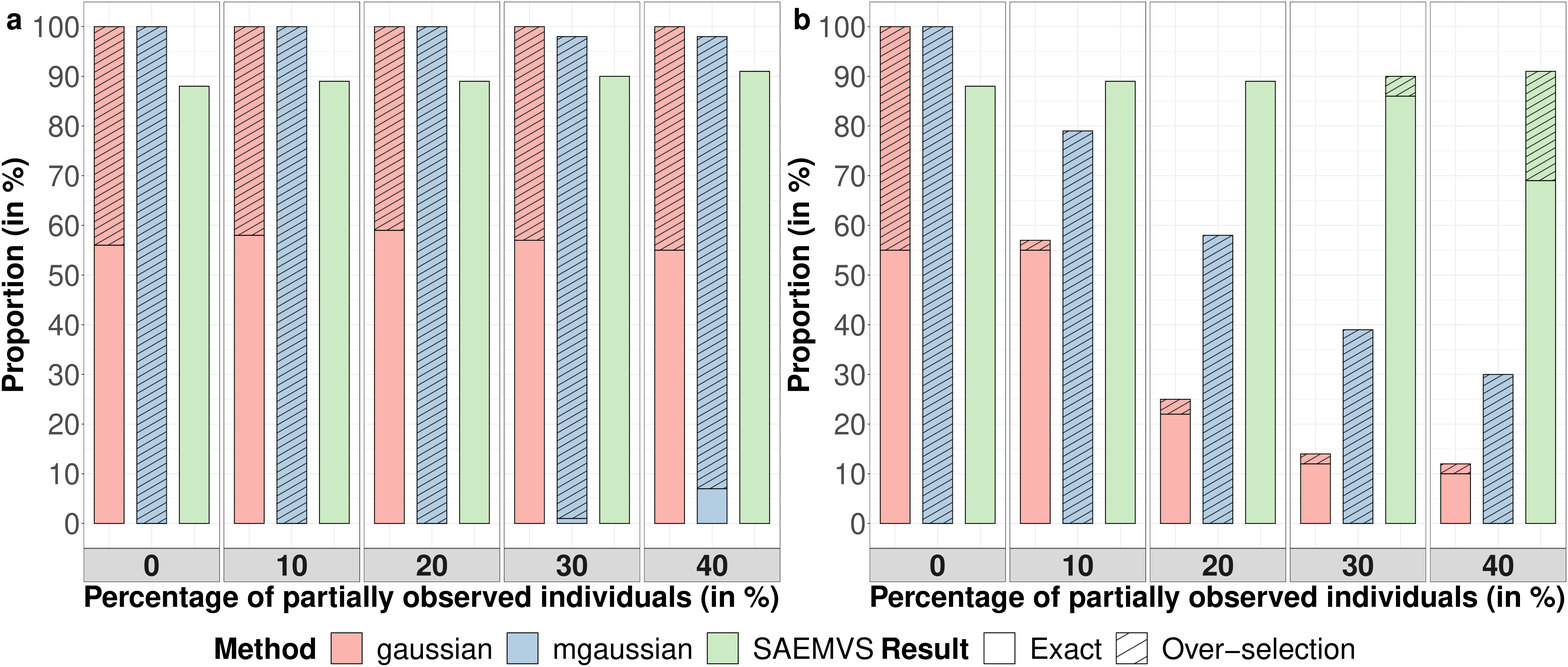}
    \caption{Proportion of data-sets on which the three methods (in colour) select the correct model ("Exact", unpatterned bars), or a model that strictly includes the correct model ("Over-selection", striped bars) for the first individual parameter (a) and the second individual parameter (b), and different percentage of partially observed individuals (on the $x$-axis).} 
    \label{Plots_phi1_phi2}
\end{figure}

Thus, it is expected that the number of partially observed individuals will have a negative impact on the estimated support of the second individual parameter. Indeed, this is what is observed in Figure \ref{Plots_phi1_phi2} for gaussian and mgaussian methods. On this figure, for the first individual parameter (graph a), we can see that gaussian and mgaussian methods select a model that almost always includes the true support (striped bars). However, mgaussian method almost never selects the true model and the gaussian method only in one case out of 2. In contrast, SAEMVS selects exactly the right model (unpatternedbars) in a large majority of cases (about $90\%$) with no false positives. Note that in many applications, especially in biology, false positives are to be avoided due to the cost of the experiments, and therefore SAEMVS seems to be more efficient in this regard. For the second individual parameter (graph b), the methods gaussian and mgaussian suffer greatly from the increase in the number of partially observed individuals, whereas the mixed effect model structure in SAEMVS shows greater robustness thanks to the classical pooling of information phenomenon, whereby individuals with missing data can benefit from the remaining fully observed individuals. 

\renewcommand\figurename{Fig.}

\subsection{Impacts of the different parameters and collinearity between covariates}
\label{simul_onedim}

\subsubsection{Model}
\label{sim_model}

Now, the performance of SAEMVS is studied under a variety of scenarios. This exploration is carried out on a different nonlinear model, highlighting our approach's flexibility. For ease of presentation, we consider variable selection in a one-dimensional setting, \textit{i.e.} $q=1$. The other parameters are assumed to be shared among individuals and they are estimated jointly. A data-set is simulated according to a logistic growth model such as: 
\begin{equation}
\label{model_simul}
    \left\{
    \begin{array}{ll}
        y_{ij} = \dfrac{\psi_1}{1+\exp\left(-\dfrac{t_{ij}-\varphi_i}{\psi_2}\right)} + \varepsilon_{ij} ,& \varepsilon_{ij} \overset{\text{i.i.d.}}{\sim} \mathcal{N}(0,\sigma^2), \\
        \varphi_i = \mu + \beta^\top V_i+ \xi_i ,& \xi_i \overset{\text{i.i.d.}}{\sim} \mathcal{N}(0,\Gamma^2),
    \end{array}
\right.
\end{equation}
where $\psi=(\psi_1,\psi_2)$ is seen as unknown fixed effects, $\varphi_i\in \mathbb{R}$, $\mu \in \mathbb{R}$, $\beta \in \mathbb{R}^p$, and $\Gamma^2>0$. Thus, parameter $\psi$ must also be estimated and therefore the population parameter is $\theta=(\mu,\beta,\psi,\Gamma^2,\sigma^2)$. The procedure presented earlier in this paper can easily be extended to this case. 

Indeed, as the function $g$ here is not separable into $\varphi_i$ and $\psi$, the model does not belong to the curved exponential family since it is not possible to write $\overset{\sim}{Q}_{1}$ with an exponential form. As a result, the expression for the maximum argument in $\psi$ at M-step is not explicit. One solution would be to do numerical optimisation in $\psi$. However, for ease of implementation of the MCMC-SAEM algorithm, following the idea of \cite{kuhn2005}, an extended model belonging to the curved exponential family is used to estimate the parameters by considering:
\begin{equation}
\label{extended_model}    \left\{
    \begin{array}{ll}
        y_{ij} &\overset{\text{ind.}}{\sim} \mathcal{N}(g(\varphi_i, \psi, t_{ij}),\sigma^2), \\
        \varphi_i &\overset{\text{ind.}}{\sim} \mathcal{N}(\mu + \beta^\top V_i,\Gamma^2), \\
        \psi &\sim \mathcal{N}(\eta,\Omega),
    \end{array}
\right.
\end{equation}
with $\varphi$ and $\psi$ independent, $\Omega=\textrm{diag}(\omega^2_1,\omega^2_2)$ known and $\theta^{ext}=(\mu,\beta,\eta,\sigma^2,\Gamma^2)$ the new population parameter to be estimated. The estimation of $\eta$ is then used as an estimation of $\psi$. As previously, the indicators $\delta=(\delta_{\ell})_{1 \leq \ell \leq p}$ (Equation~\eqref{delta} with $q=1$) are introduced, and we consider the same priors as in \eqref{priors} for $(\mu,\beta,\sigma^2,\Gamma^2,\delta,\alpha)$ but in their one-dimensional version: 
\begin{equation}
\label{priors_unidim}
\begin{array}{ll}
    \pi(\beta \vert \delta) &= \mathcal{N}_{p}(0,D_{\delta})\text{, with } D_{\delta}=\textrm{diag}((\mathbf{1}-\delta)\nu_0 + \delta \nu_1) \text{, } 0 < \nu_0 < \nu_1 \text{, } \\
    \pi(\mu)&=  \mathcal{N}(0,\sigma^2_{\mu}) \text{, with } \sigma^2_{\mu}>0 \text{, } \\
    \pi(\sigma^2)&= \mathcal{IG}\left(\dfrac{\nu_{\sigma}}{2},\dfrac{\nu_{\sigma}\lambda_{\sigma}}{2}\right) \text{, with } \nu_{\sigma},\lambda_{\sigma}>0 \text{, } \\
    \pi(\Gamma^2)&= \mathcal{IG}\left(\dfrac{\nu_{\Gamma}}{2},\dfrac{\nu_{\Gamma}\lambda_{\Gamma}}{2}\right) \text{, with } \nu_{\Gamma},\lambda_{\Gamma}>0 \text{, } \\
    \pi(\delta\vert\alpha)&= \alpha^{\vert\delta\vert}(1-\alpha)^{p-\vert\delta\vert} \text{, with } \alpha \in [0,1] \text{ and } \vert\delta\vert=\sum_{\ell=1}^{p} \delta_{\ell} \text{, } \\
    \pi(\alpha) &= \textit{Beta}(a,b) \text{, with } a,b>0. 
\end{array}
\end{equation}
For $\eta$, the following prior is chosen: for $r \in \{1, 2\}$, $\pi(\eta_r) = \mathcal{N}(0,\rho_r^2)$, with $\rho_r^2>0$ known. This amounts to randomising hyperparameters of the prior on $\psi$, implying a less informative prior than if $\eta$ were fixed. Here, $\Theta=(\theta^{ext}, \alpha)$ is the population parameter and $Z=(\varphi,\psi,\delta)$ are the latent variables. The steps of the MCMC-SAEM algorithm can be adapted to this model. The estimation method is unchanged for parameters $(\mu,\beta,\sigma^2,\Gamma^2,\alpha)$ because the quantity $\overset{\sim}{Q}_{1}$ is separable into $(\mu,\beta,\sigma^2,\Gamma^2,\alpha)$ and $\eta$. The main difference is that there is another latent variable $\psi$, which must also be simulated at the S-step. The thresholding procedure is not modified. See Appendix \ref{app_extension_fixed} for more details about this extension of SAEMVS.

\begin{remark}
\label{rmk_omega}
In order to limit the estimation error between the initial model and this extended model, the value of the covariance matrix is adapted during the iterations. Inspired by the results of \cite{allassonniere_curved_2021} for the case of the computation of the MLE, the following process is chosen: start with a fairly large initial value $\Omega^{(0)}~=~\textrm{diag}(\omega_1^{2^{(0)}},\omega_2^{2^{(0)}})$ for a certain number $\kappa$ of iterations, then multiply it by $0<\tau<1$, and iterate this process every $\kappa$ iterations. Starting from a large initial value, the value of $\Omega$ remains large enough during the first iterations to allow a rather fast convergence speed, then it is slowly decreased towards $0$, while remaining always strictly positive, to limit the estimation error between the initial model and the extended model.
\end{remark}

\subsubsection{Simulation design}
\label{sim_settings}

For this simulation study, individual profiles are simulated according to model~\eqref{model_simul} by considering $n_1 = \dots = n_n =10$ observations per individual and regular observation time points such that $t_{ij}=t_j=150+(j-1)\dfrac{3000-150}{J-1}$, $\sigma^2=30$, $\psi_1=200$, $\psi_2~=~300$, $\mu=1200$, $\beta=(100,50,20,0, \dots, 0)^\top$. 
Thus, only the first three covariates are assumed to be influential and their respective intensities are contrasted. The individual covariates $V_i \in \mathbb{R}^p$, $1 \leq i \leq n$, are simulated independently according to a centred multivariate Gaussian distribution with covariance matrix $\Sigma \in \mathcal{M}_p(\mathbb{R})$. To test the sensitivity of SAEMVS to the correlation that may exist between covariates, different scenarios are tested corresponding to different structures for matrix $\Sigma$. Different values of $n$ (number of subjects) and $p$ (number of covariates) are used according to the scenario. 
Several values of $\Gamma^2$ (variance of the random effects) are also used in order to evaluate the performances of SAEMVS in different "signal-to-noise ratio" situations.

\begin{itemize}
\item \textbf{Scenario with uncorrelated covariates.} This is the baseline scenario where optimal performance of Algorithm~\ref{algoSAEMVS} is expected. This corresponds to $\Sigma=I_p$, where $I_p$ is the identity matrix of size $p$.  The following values for $n$, $p$ and $\Gamma^2$ are used: $n \in \{100,200\}$, $p \in \{500,2000,5000\}$ and $\Gamma^2 \in \{200, 1000, 2000 \}$.

\item \textbf{Scenarios with correlations between covariates.}

\begin{enumerate}
    \item The first scenario leaves the three influential covariates uncorrelated with all other covariates whereas the non-influential covariates are correlated with each other. An autoregressive correlation structure is considered between the non-influential covariates. This corresponds to $\Sigma = \left (
   \begin{array}{c|c}
      I_3 & 0_{3,p-3} \\
      \hline
      0_{p-3,3} & (\rho_{\Sigma}^{\vert i-j\vert})_{i,j \in \{4,\dots,p\}} \\ 
   \end{array}
\right)$, with $\vert\rho_{\Sigma}\vert<1$. 
    \item In the second scenario, the third influential covariate is assumed to be correlated to every non-influential covariate according to an autoregressive correlation structure. This corresponds to 
    $\Sigma = \left (
   \begin{array}{c|c}
      I_3 & A \\
      \hline
      A^\top & I_{p-3} \\ 
   \end{array}
\right)$, with\\ $A=\left (
   \begin{array}{cccc}
      0 & 0 & \dots & 0 \\
      0 & 0 & \dots & 0 \\
      & & (\rho_{\Sigma}^{\vert 3-j\vert})_{j \in \{4,\dots,p\}} & \\ 
   \end{array}
\right)$, $\vert\rho_{\Sigma}\vert<1$. 
    \item The third scenario considers correlations between the sole influential covariates. Again, an autoregressive correlation structure is used. This corresponds to $\Sigma~=~\left (
   \begin{array}{c|c}
      (\rho_{\Sigma}^{\vert i-j\vert})_{i,j \in \{1,\dots,3\}} & 0_{3,p-3} \\
      \hline
      0_{p-3,3} & I_{p-3} \\ 
   \end{array}
\right)$, $\vert\rho_{\Sigma}\vert<1$. 
    \item In the fourth scenario, an autoregressive correlation structure is used between the covariates without making any distinction between the influential covariates and the non-influential covariates. This corresponds to $\Sigma~=~(\rho_{\Sigma}^{\vert i-j\vert})_{i,j \in \{1,\dots,p\}}$, $\vert\rho_{\Sigma}\vert<1$.
\end{enumerate}
\end{itemize}

To study the impact of correlations between covariates according to the four scenarios above, the following values for $n$, $p$, $\Gamma^2$ and $\rho_{\Sigma}$ are used: $n=200$, $p \in \{500,2000,5000\}$, $\Gamma^2 \in \{200,2000\}$ and $\rho_{\Sigma} \in \{0.3,0.6\}$. 

For each of the five scenarios described above and each combination $(n,p,\Gamma^2)$ or $(n,p,\Gamma^2,\rho_{\Sigma})$, $100$ different data-sets are simulated and the support of $\beta$ is estimated by applying Algorithm~\ref{algoSAEMVS} on each data-set. The algorithmic settings are given in Appendix \ref{algo_setting}. Note that, in order to be able to compare covariates that do not have the same order of magnitude, the covariates are standardised. The performances in terms of exact selection of the true influential covariates, over-selection and under-selection are examined (see Subsection~\ref{sim_results}). 

\subsubsection{Results}
\label{sim_results}

\begin{figure}[!ht]
    \centering
    \includegraphics[width=0.88\textwidth]{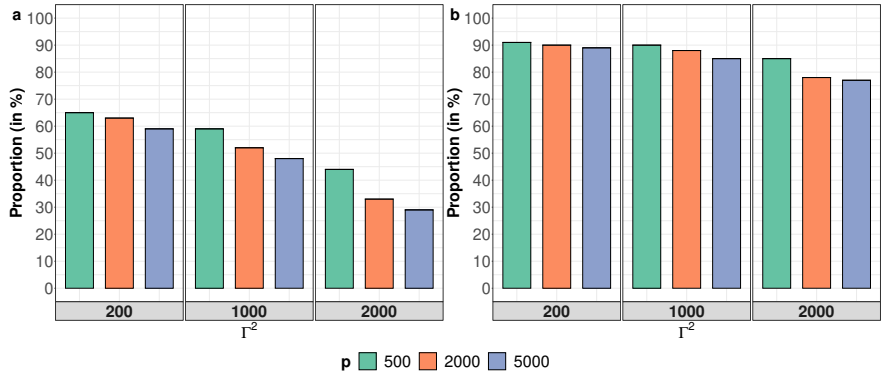}
    \caption{Uncorrelated covariates. Proportion of data-sets on which SAEMVS selects the correct model for $n=100$ (a) and $n=200$ (b), and different values of $p$ and $\Gamma^2$.} 
    \label{exact_model}
\end{figure}

\paragraph{Scenario with uncorrelated covariates}
The results are presented in Table~\ref{Tab_ind} and Figure~\ref{exact_model}. SAEMVS selects exactly the right model in a large majority of cases for a sufficiently large number of individuals $n$. When $n$ increases, the results improve, which suggests a consistency property in selection. With $n$ and $p$ fixed, the more the inter-individual variance $\Gamma^2$ is important, the more the results degrade. Indeed, as $\Gamma^2$ increases, the "signal-to-variability" ratio decreases, leading to 
difficulties in detecting the third covariate associated with the lowest non-zero coefficient in $\beta$. 
It could also be noted that with $n$ and $\Gamma^2$ fixed, the results deteriorate when $p$ increases but the effect of $p$ seems weak when $n$ is large. 
In addition, when SAEMVS fails, it is most often because it under-selects, that is, it selects fewer variables than there are. Indeed, average sensitivity values are lower than those for specificity. It seems that SAEMVS tends to avoid false positives, even though this may result in not having selected all the truly influential covariates. 

\begin{table*}[!ht] 
\begin{minipage}{\textwidth} 
\centering
\caption{Uncorrelated covariates. For each of the quantifiers (Sensitivity, Specificity, Accuracy), the mean over the $100$ data sets is shown, with the empirical standard deviation divided by $\sqrt{100}$ in brackets.}\label{Tab_ind}%
\scriptsize
\def\thefootnote{\alph{footnote}}
\begin{tabular}{@{}lllllllll@{}}
\toprule
& & \multicolumn{3}{c}{$n=100$} & &\multicolumn{3}{c}{$n=200$}\\
\cline{3-5} 
\cline{7-9}
$\Gamma^2$ & $p$ & Se\footnotemark[1] & Sp\footnotemark[2] & Ac\footnotemark[3] & & Se\footnotemark[1] & Sp\footnotemark[2] & Ac\footnotemark[3] \\
\midrule
200 & 500 & 0.883 & 1 & 0.999 & & 0.973 & 1 & 1 \\
     &     & (0.0180) & {(5e-05)} & {(0.0001)} & & (0.0091) & {(2e-05)} & {(6e-05)} \\ 
\addlinespace[2ex]
200 & 2000 & 0.900 & 1 & 1 & & 0.987 & 1 & 1 \\
     &      & (0.0160) & {(2e-05)} & {(3e-05)} & & (0.0066) & {(2e-05)} & {(2e-05)} \\
\addlinespace[2ex]
200 & 5000 & 0.870 & 1 & 1 & & 0.983 & 1 & 1 \\
     &      & (0.0189) & {(6e-06)} & {(1e-05)} & & (0.0073) & {(1e-05)} & {(1e-05)} \\
\midrule
1000 & 500 & 0.860 & 1 & 0.999 & & 0.977 & 1 & 1 \\
      &      & (0.0185) & {(3e-05)} & {(0.0001)} & & (0.0086) & {(3e-05)} & {(6e-05)} \\ 
\addlinespace[2ex]
1000 & 2000 & 0.840 & 1 & 1 & & 0.977 & 1 & 1 \\
      &      & (0.0180) & {(9e-06)} & {(3e-05)} & & (0.0086) & {(1e-05)} & {(2e-05)} \\
\addlinespace[2ex]
1000 & 5000 & 0.813 & 1 & 1 & & 0.963 & 1 & 1 \\
      &      & (0.0197) & {(4e-06)} & {(1e-05)} & & (0.0105) & {(6e-06)} & {(8e-06)} \\
\midrule
2000 & 500 & 0.810 & 1 & 0.999 & & 0.950 & 1 & 1 \\
      &      & (0.0185) & {(3e-05)} & {(0.0001)} & & (0.0120) & {(0)} & {(7e-05)} \\ 
\addlinespace[2ex]
2000 & 2000 & 0.777 & 1 & 1 & & 0.937 & 1 & 1 \\
      &      & (0.0178) & {(2e-05)} & {(3e-05)} & & (0.0131) & {(1e-05)} & {(2e-05)} \\
\addlinespace[2ex]
2000 & 5000 & 0.767 & 1 & 1 & & 0.930 & 1 & 1 \\
      &      & (0.0174) & {(6e-06)} & {(1e-05)} & & (0.0137) & {(4e-06)} & {(9e-06)} \\
\bottomrule
\end{tabular}
\footnotetext[1]{Sensitivity is the proportion of true positives correctly identified.}
\footnotetext[2]{Specificity is the proportion of the true negatives correctly identified.}
\footnotetext[3]{Accuracy is the proportion of true results, either true positive or true negative.}
\end{minipage}
\end{table*}

\paragraph{Scenarios with correlated covariates}

\begin{figure}[!ht]
  \begin{center}
  \includegraphics[width=0.88\textwidth]{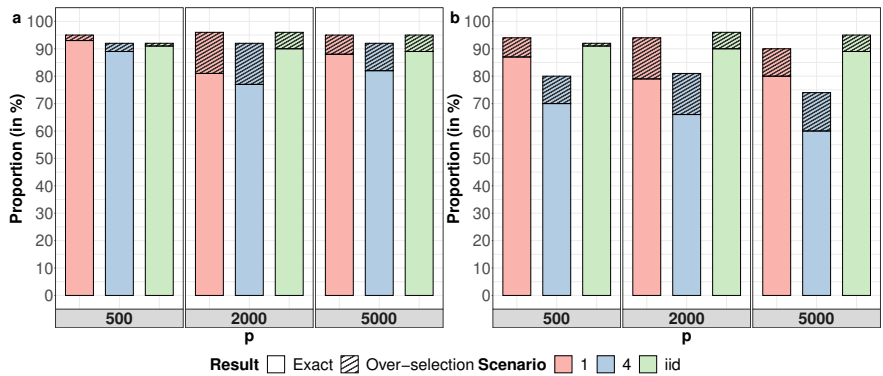}
    \caption{Correlated covariates for $\Gamma^2=200$. Proportion of data-sets on which SAEMVS selects the correct model ("Exact", unpatterned bars) or a model that strictly includes the correct model ("Over-selection", striped bars) for $\rho_{\Sigma}=0.3$ (a), $\rho_{\Sigma}=0.6$ (b), and different values of $p$. Scenario "i.i.d" corresponds to the case where the covariates are not correlated and is used as a reference.} 
    \label{comp_corr}
  \end{center}
\end{figure}

The results for scenarios $1$ and $4$ with $\Gamma^2=200$ are presented in Figure~\ref{comp_corr} and Table~\ref{Tab_corr} for the two values of $\rho_{\Sigma}$. The results for the other scenarios are given in Appendix \ref{ann_corr}. On these figures, one can compare the selection performance of SAEMVS in the different scenarios of correlations between covariates with the case without correlations (i.i.d scenario). First, for scenario $1$, that is when the non-active covariates are correlated, quite similar performances to the i.i.d scenario are observed, but with more over-selection. Indeed, as a consequence of the correlation between irrelevant covariates, the latter tend to be selected more often and in small groups. Then, scenario $4$ corresponds to a full correlation matrix between all covariates. Note that the correlation matrix chosen for this scenario assumes a fairly strong correlation between the three true covariates. Thus, as these active covariates explain the response variable in a similar way, this scenario leads to much under-selection compared to the i.i.d case. This can be seen in the sensitivity values. Moreover, it over-selects more than scenario i.i.d because of the correlations between the true and false covariates. However, by looking at the specificity values, we can see that even in cases where SAEMVS over-selects, the false positive rate (which is equivalent to (1 – specificity)) remains very close to $0$. 

\begin{table*}[!ht] 
\begin{minipage}{\textwidth} 
\centering
\caption{Correlated covariates for $\Gamma^2=200$. For each of the quantifiers (Sensitivity, Specificity, Accuracy), the mean over the $100$ data sets is shown, with the empirical standard deviation divided by $\sqrt{100}$ in brackets.}\label{Tab_corr}%
\scriptsize
\def\thefootnote{\alph{footnote}}
\begin{tabular}{@{}lllllllll@{}}
\toprule
& & \multicolumn{3}{c}{$\rho_{\Sigma}=0.3$}&  &\multicolumn{3}{c}{$\rho_{\Sigma}=0.6$}\\
\cline{3-5} 
\cline{7-9}
Scenario & $p$ & Se\footnotemark[1] & Sp\footnotemark[2] & Ac\footnotemark[3] & & Se\footnotemark[1] & Sp\footnotemark[2] & Ac\footnotemark[3] \\
\midrule
i.i.d & 500 & 0.973 & 1 & 1 & & 0.973 & 1 & 1 \\
   &  & (0.0091) & {(2e-05)} & {(6e-05)}  & & (0.0091) & {(2e-05)} & {(6e-05)} \\ 
\addlinespace[2ex]
1 & 500 & 0.983 & 1 & 1 & & 0.980 & 1 & 1 \\
  &   & (0.0073) & {(3e-05)} & {(5e-05)} & & (0.0080) & {(6e-05)} & {(8e-05)} \\ 
\addlinespace[2ex]
4 & 500 & 0.973 & 1 & 1 & & 0.933 & 1 & 0.999 \\
  &   & (0.0091) & {(4e-05)} & {(7e-05)} & &  (0.0134) & {(0.0001)} & {(0.0002)} \\ 
\addlinespace[2ex]
\midrule
i.i.d & 2000 & 0.987 & 1 & 1 & & 0.987 & 1 & 1 \\
  &   & (0.0066) & {(2e-05)} & {(2e-05)} & & (0.0066) & {(2e-05)} & {(2e-05)} \\ 
\addlinespace[2ex]
1 & 2000 & 0.987 & 1 & 1 & & 0.980 & 1 & 1\\
  &   & (0.0066) & {(3e-05)} & {(3e-05)} & & (0.0080) & {(3e-05)} & {(3e-05)}\\ 
\addlinespace[2ex]
4 & 2000 & 0.973 & 1 & 1 & & 0.937 & 1 & 1\\
  &   & (0.0091) & {(2e-05)} & {(3e-05)} & & (0.0131) & {(3e-05)} & {(4e-05)} \\ 
\addlinespace[2ex]
\midrule
i.i.d & 5000 & 0.983 & 1 & 1 & & 0.983 & 1 & 1 \\
  &   & (0.0073) & {(1e-05)} & {(1e-05)} & & (0.0073) & {(1e-05)} & {(1e-05)}\\ 
\addlinespace[2ex]
1 & 5000 & 0.983 & 1 & 1 & & 0.967 & 1 & 1\\
  &   & (0.0073) & {(1e-05)} & {(1e-05)} & & (0.0101) & {(9e-06)} & {(1e-05)} \\ 
\addlinespace[2ex]
4 & 5000 & 0.973 & 1 & 1 & & 0.913 & 1 & 1  \\
   &  & (0.0091) & {(8e-06)} & {(9e-06)} &  & (0.0147) & {(1e-05)} & {(2e-05)}\\ 
\addlinespace[2ex]
\bottomrule
\end{tabular}
\footnotetext[1]{Sensitivity is the proportion of true positives correctly identified.}
\footnotetext[2]{Specificity is the proportion of the true negatives correctly identified.}
\footnotetext[3]{Accuracy is the proportion of true results, either true positive or true negative.}
\end{minipage}
\end{table*}

\subsection{Comparison with an MCMC implementation}
\label{comp_MCMC}

It is reasonably straightforward to implement an MCMC algorithm for full posterior inference on the spike-and-slab variable selection for non-linear mixed-effects model.
To understand precisely the added value of the SAEM algorithm, we
compare the run time of a full MCMC approach and the MCMC-SAEM method proposed in this paper, and highlight the better scaling properties of the latter. 
To build the most informative comparison, the same model with a smooth spike is considered for both the MCMC and MCMC-SAEM approaches, remarking that spike-and-slab priors with a Dirac spike are known to pose challenges for MCMC \citep[see][]{baiSpikeandSlabMeetsLasso2021a}. 
For the MCMC algorithm, an efficient C++ implementation of a random walk adaptive MCMC is used through the Nimble software \citep{devalpineProgrammingModelsWriting2017}, which uses an adaptive scheme proposed in \cite{shaby2010exploring}.
To make the comparison as fair as possible, we marginalise the sampler over the discrete inclusion variables $\delta$, to mirror the marginalisation in \eqref{Qtilde}. 
This was found to appreciably improve the mixing of the MCMC algorithm. It is possible to retrieve the $\delta$ variables from the posterior samples using their conditional posterior distribution.

Common data-sets are simulated according to model~\eqref{model_simul} with the following parameters: $n=200$ individuals, $n_1~=~\dots~=~n_n~=~10$ observations per individual, $p~\in~\{500, 700, 1000, 1500, 2000, 2500 \}$ covariates, $\sigma^2=30$, $\psi_1=200$, $\psi_2~=~300$, $\mu=1200$, $\beta~=~(100,50,20,0, \dots, 0)^\top$ and $\Gamma^2=200$. For $ i \in \{1,\dots ,n\}$ and $j \in \{1,\dots,10\}$, $t_{ij}~=~t_j=150+(j-1)\dfrac{3000-150}{J-1}$. Covariates are simulated independently and identically distributed according to $\mathcal{N}(0,1)$. 
The objective is to compare the time needed to estimate the parameters $\Theta~=~(\mu,\beta,\psi,\Gamma^2,\sigma^2,\alpha)$ between the MCMC-SAEM algorithm proposed in this article (Algorithm~\ref{algoMCMCSAEM} adapted to model~\eqref{model_simul}, see Appendix \ref{app_extension_fixed}) and the full MCMC procedure described above. As explained in Subsection~\ref{sim_model}, to estimate the parameters, we consider the extended model~\eqref{extended_model}. 
The same model structure~\eqref{extended_model} and priors are used for both approaches. For $(\mu,\beta,\Gamma^2,\delta,\alpha)$ the priors are as in \eqref{priors_unidim}. For $\eta$, the prior of Subsection~\ref{sim_model} is chosen: for $r \in \{1, 2\}$, $\pi(\eta_r) = \mathcal{N}(0,\rho_r^2)$, with $\rho_r^2>0$ known. To stabilise the MCMC procedure, the prior on $\sigma^2$ is modified to a uniform distribution on $[0,200]$ for both methods. This has very little consequence for SAEMVS. Indeed, the only difference lies in the updating of $\sigma^2$ at the M-step of the MCMC-SAEM algorithm, which becomes: 
$$\sigma^{2^{(k+1)}}=\left\{
    \begin{array}{ll}
        \dfrac{s_{1,k+1}}{n_{tot}} & \mbox{if } \dfrac{s_{1,k+1}}{n_{tot}}\leq 200 \\
        200 & \mbox{else.}
    \end{array}
\right. $$
The two methods are both initialised with: $\forall \ell \in \{1,\dots, 10\}$, $\beta^{(0)}_{\ell}=100$, $\forall \ell~\in~\{11,\dots, p\}$, $\beta^{(0)}_{\ell}=1$, $\mu^{(0)}=1400$, $\sigma^{2^{(0)}}=100$, $\alpha^{(0)}=0.1$ and $\eta^{(0)}~=~(400,400)^\top$. In practice, to avoid convergence toward a local maximum in the MCMC-SAEM algorithm, a simulated annealing version of SAEM \citep[see][]{lavielle} is implemented. Thus, in SAEMVS, $\Gamma^2$ is initialised very large to explore the space during the first iterations, with $\Gamma^{2^{(0)}}=5000$. For the full MCMC procedure, a more plausible value of $\Gamma^2$, $\Gamma^{2^{(0)}}=500$, is chosen as initialisation. The hyperparameters are set in the same way for both methods as well: $\nu_0=0.04$, $\nu_1=12000$, $\sigma_{\mu}=3000$, $\nu_{\Gamma}=\lambda_{\Gamma}=1$, $a=1$, $b=p$, $\Omega=\textrm{diag}(20,20)$ and $\rho_1^2=\rho_2^2=1200$. 

\begin{figure}[ht]
    \centering
    \includegraphics[width=0.8\textwidth]{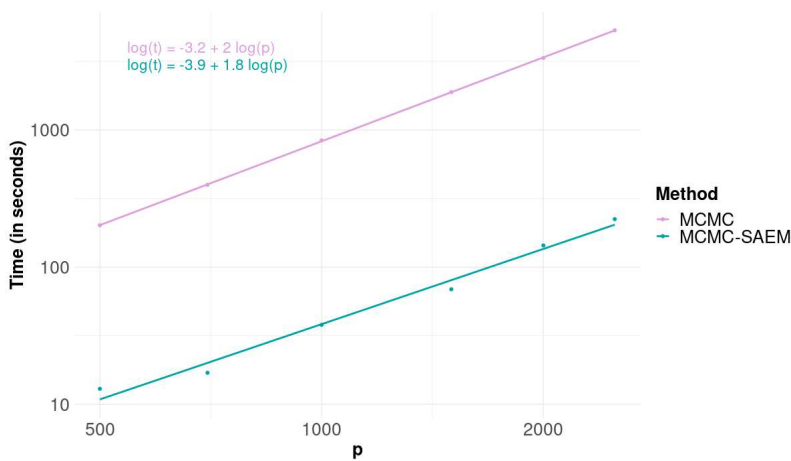}
    \caption{Comparison of computation times between MCMC (in purple) and Algorithm~\ref{algoMCMCSAEM} (MCMC-SAEM, in blue) inference methods in $\log_{10}$-$\log_{10}$ scale.}
    \label{fig:comp_temps}
\end{figure}

For ease of presentation, we compare the two approaches for a single value of $\nu_0$. It is standard practice to run an MCMC spike-and-slab model for a single value (see for instance \cite{george97} or \cite{malsiner2018comparing}). 
The MCMC algorithm was run for $3000$ iterations, which was just enough to reach convergence (assessed by comparing multiple chains) for a variety of $\nu_0$ and $p$ values. The MCMC-SAEM algorithm was run for $500$ iterations and showed appropriate convergence. This convergence was assessed qualitatively by looking at the difference between successive parameter estimates during iterations. Under these conditions, for all $p~\in~\{500, 700, 1000, 1500, 2000, 2500 \}$, both methods were run for $50$ different data-sets and the minimum time was kept for each method. The results obtained are shown in Figure~\ref{fig:comp_temps}. In this figure, computation times of the full MCMC procedure (in purple) and of MCMC-SAEM (in blue) are represented by the points for the different values of $p$. The lines represent the regression line associated with each method. Note that a $\log_{10}$-$\log_{10}$ scale is used in this figure. This shows that both methods have an execution time that grows polynomially with $p$. Furthermore, the polynomial complexity of the two methods, \textit{i.e.} the slope of the regression lines, is slightly lower for the MCMC-SAEM method. Thus, if we note respectively $\tau_{MCMC}$ and $\tau_{MCMC-SAEM}$ the execution time associated with each of the methods under the conditions previously described, empirically Figure~\ref{fig:comp_temps} strongly suggests that $\dfrac{\tau_{MCMC}}{\tau_{MCMC-SAEM}}\approx 10^{0.7} p^{0.2}$. To sum up, the MCMC-SAEM algorithm proposed in this paper appears $10^{0.7} p^{0.2}$ times faster than the classical MCMC procedure, \textit{i.e.} between $17$ and $24$ times faster for $p$ between $500$ and $2500$. In other words, SAEMVS allows to browse a grid of about $20$ values of the penalisation parameter $\nu_0$ while a classical MCMC only looked at one value of this parameter.

\section{Application to plant senescence genetic marker identification}
\label{Section_data_reel}

In this section, SAEMVS is applied to a problem of marker-assisted selection to assist breeding of winter wheat. 
We are interested in identifying genetic markers that impact the senescence process, \textit{i.e.} the ageing of these plants, which is under genetic determinism.
Because heading date is related to senescence, the relevance of the selected senescence markers is then compared to known flowering genes as well as to markers associated to heading data obtained by applying an association mapping model \citep{yu} to these data. In the following, these markers are called heading QTLs. Note that the heading date can be seen as an easily measurable approximation of the flowering date.

\subsection{Data description and pre-processing}

The plant material used here has been previously described in \cite{rincent2018,rincent2019} and \cite{touzy}. It is composed of $n=220$ wheat varieties. Senescence was measured as the global proportion of senesced surface of the canopy. This proportion was observed on each variety $J=18$ times over time. 

For each of the $n$ varieties, information from high-throughput genotyping is available on several tens of thousands of SNPs positioned along the entire genome \citep{rimbert} (see Appendix \ref{app_genotyping}). 
These binary variables constitute the covariates to be selected in order to explain the differences between varieties in terms of senescence. 
There are structurally strong correlation and collinearity between these variables, which hampers variable selection, as pointed out in Section~\ref{simul_onedim} and discussed in many references (see e.g. \cite{malsiner2018comparing} and \cite{heuclin2021bayesian} in spike and slab models). 
To permit a comparison with the genetic markers associated with flowering in the context of collinearity, SAEMVS is applied chromosome by chromosome. 
Note, however, that variable selection on the whole genome (~$\approx$~30000 markers) would be computationally feasible in a matter of hours with strict pre-processing to address multicollinearity/near-multicollinearity. 
We still apply the following minimal pre-processing: if several SNPs in the same chromosome are strictly collinear for all varieties, only one of them is retained for analysis. 
This eliminates $916$ SNPs across the genome. 
As a result, for a given chromosome $C$, the number $p$ of covariates in $V_i^C$ is always less than $2000$, while remaining in a high-dimensional framework where $p >> n$ (see Table~\ref{Tab_chr} in Appendix \ref{annexe_table} where the values of $p$ are given for each chromosome).
The wheat varieties in this data-set are structured into genetic groups, which may cause confusion between the effect of the subpopulation structure and that of SNPs. 
To control for subpopulation structure, we adapt our model to include covariates not subject to selection. We consider the first 5 principal components of a principal coordinates analysis performed on the available SNPs in the model. 
These 5 adjustment variables, capturing subpopulation structure, are denoted by $v_i\in \mathbb{R}^5$ for variety $i$ and are guaranteed to be used in the regression.

\subsection{Modelling}
\label{modeling_data_reel}

Some senescence curves are shown in Figure~\ref{data_reel} in Appendix \ref{app_fig_data_reel}. We can see that a logistic growth model (Equation~\eqref{model_simul}) with a maximum value of 100\% and inter-individual variability on the two other parameters is coherent with the shape of these curves. This is a simple model for the purposes of the example, but note that the variable selection method can accommodate other models like beta regression models with minimal changes. Here, we choose to analyse the effects of SNPs from each chromosome $C$ on a single parameter of interest: the variability of characteristic times between varieties.
Denoting $y_{ij}$ the proportion of senesced surface of the plant of the variety $1\leq i \leq n$ at time $t_{ij}$, this leads to the following model when focussing on chromosome C: 
\begin{equation*}
\label{model_data_reel}
    \left\{
    \begin{array}{ll}
        y_{ij} = \dfrac{100}{1 + \exp\left(-\dfrac{t_{ij}-\varphi_i}{\psi_i}\right)} + \varepsilon_{ij} ,& \varepsilon_{ij} \overset{\text{i.i.d.}}{\sim} \mathcal{N}(0,\sigma^2), \\
        \varphi_i = \mu + \lambda^\top v_i + \beta^\top V_i^C+ \xi_i ,& \xi_i \overset{\text{i.i.d.}}{\sim} \mathcal{N}(0,\Gamma^2),  \\
        \psi_i = \eta + \omega_i ,& \omega_i \overset{\text{i.i.d.}}{\sim} \mathcal{N}(0,\Omega^2).
    \end{array}
\right.
\end{equation*}
The parameter to be estimated is therefore: $\theta=(\mu,\lambda,\beta,\eta,\sigma^2,\Gamma^2,\Omega^2)$. As previously, the indicators $\delta$ (Equation~\eqref{delta} with $q=1$) are introduced and we force the inclusion of variables $v_i$ in the model by a useful reformulation similar to what was done for the intercept in \eqref{reformulation}.
We use the following priors for $\eta$ and $\Omega^2$ : $\pi(\eta) = \mathcal{N}(0,\sigma_{\eta}^2)$, with $\sigma_{\eta}^2$ known, and $\pi(\Omega^2)= \mathcal{IG}\left(\dfrac{\nu_{\Omega}}{2},\dfrac{\nu_{\Omega}\lambda_{\Omega}}{2}\right)$, with $\nu_{\Omega},\lambda_{\Omega}>0$. The same priors as in \eqref{priors_unidim} are used for the other parameters. Here, $\Theta=(\theta, \alpha)$ is the population parameter and $Z=(\varphi,\psi,\delta)$ are the latent variables, where $\varphi=(\varphi_i)_{1\leq i \leq n}$ and $\psi=(\psi_i)_{1\leq i \leq n}$. The SAEMVS procedure can be easily implemented with minor modifications to the version detailed in Appendix \ref{app_extension_fixed}. The algorithmic settings are provided in Appendix \ref{ann_settings_data_reel}.

\subsection{Results and discussion}

\begin{figure}[ht]
  \begin{center}
    \includegraphics[width=0.9\textwidth]{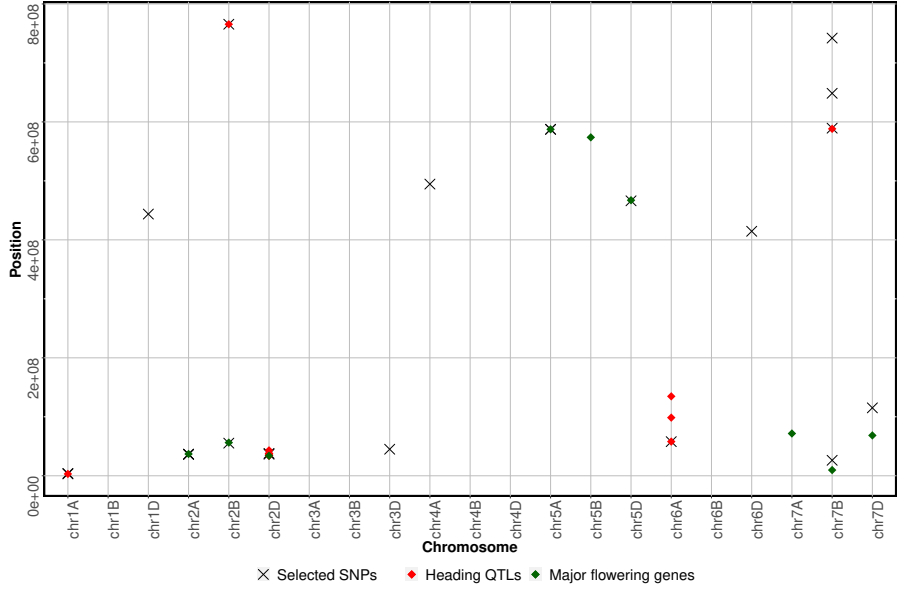}
    \caption{Position on each chromosome of the markers selected by SAEMVS (in black cross), compared to heading QTLs (in red diamond) and major flowering genes (in green diamond).}
    \label{Result_chr_by_chr}
  \end{center}
\end{figure}

The results are shown in Figure~\ref{Result_chr_by_chr} and further detailed in Table~\ref{Tab_chr} in Appendix \ref{annexe_table}. For each chromosome, the SNPs selected by SAEMVS are compared to heading QTLs and major flowering genes that were identified in previous analysis (called Ppd and Vrn genes). Indeed, as mentioned above, as flowering and senescence are linked biological processes, we can expect to find SNPs that are close to these specific QTLs and genes on the genome. As an example, SAEMVS applied to chromosome 1A leads to the selection of two SNPs that are close to each other  (\textit{i.e.} within one mega-base of each other) but also close to a heading QTL on the chromosome. On a genome-wide basis, Figure~\ref{Result_chr_by_chr} displays many colocalisations between flowering genes or heading QTLs and the SNPs selected by SAEMVS on the senescence data. Our procedure thus returns consistent results from a biological point of view. Conversely, SAEMVS also selects SNPs that are not close to zones associated with flowering on the genome, for example on chromosome 1D, 3D, 4A, 6D, 7B and 7D. The selected markers reveal potentially "stay-green" SNPs associated with late senescence independently of flowering time, which would be very interesting for plant breeders. 

\begin{figure}[ht]
    \centering
    \includegraphics[width=0.8\textwidth]{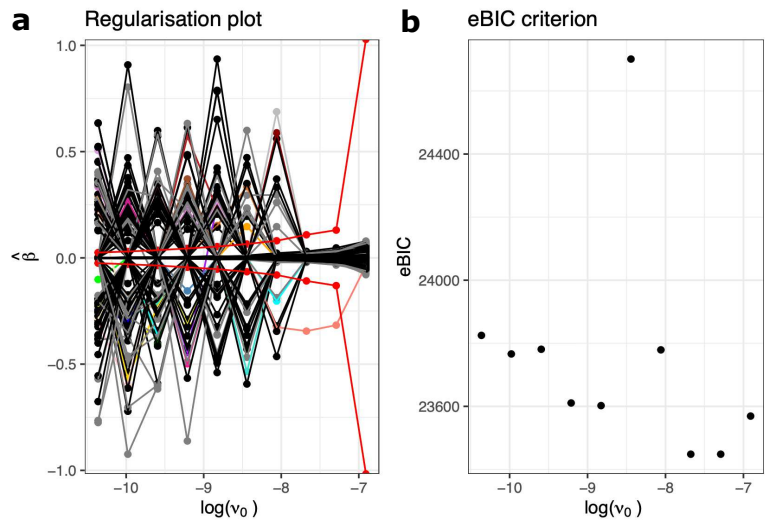}
    \caption{Regularisation plot and eBIC criterion for chromosome 6A}
    \label{chemin_regu_data_reel}
\end{figure}

It is important to note that SAEMVS suffers from selection switch among markers that are highly correlated. This is illustrated by the chevron structure in Figure~\ref{chemin_regu_data_reel}. On the left-hand plot, the red curves correspond to the selection threshold and the markers that are selected at least twice along the grid $\Delta$ of $\nu_0$ values are denoted in different colours. In such situations, we expect any variable selection to face the same challenges, as one cannot hope to select markers truly associated with the phenomenon of interest. The usual approach in biology is to rather identify regions of the genome containing markers associated to the phenomenon of interest, usually via a post-processing of the variable selection results. 

\section{Conclusion and perspectives}
\label{conclusion}

The main objective of this paper was to propose a new procedure for high-dimensional variable selection in non-linear mixed-effects models. In this work, variable selection was approached from a Bayesian perspective and a selection procedure combining the use of spike-and-slab Gaussian mixture prior and the SAEM algorithm was proposed. 
The spike-and-slab prior on the regression coefficients allows both the shrinkage towards zero of small non-significant coefficients through the spike distribution, while the largely uninformative slab distribution allows estimating influential covariates without bias from the penalisation. 
The speed of the SAEM algorithm allows exploring different levels of sparsity in the model through the variance of the spike distribution $\nu_0$, which we observed to be beneficial in selecting sparse models for regression. 
Varying the level of sparsity provides a collection of good models among which we select the one minimising an eBIC criterion.

The SAEMVS method can do both one-dimensional and multidimensional variable selection, and is does not assume the same support for each dimension. It is very flexible and was illustrated on three different models. The proposed methodology showed very good selection performance on simulated data. Indeed, SAEMVS appears to select the right support in a large majority of cases. As expected, for different numbers of covariates $p$ fixed, the right support is selected more often as the number of individuals $n$ increases and the inter-individual variance $\Gamma^2$ decreases. Even more interesting, this method is much faster than an MCMC stochastic search alternative and can solve higher-dimensional variable selection problems. The application of SAEMVS on a real data-set shows, on the one hand, the flexibility of the procedure, and on the other hand, convincing results from a biological point of view despite strong correlations/multicollinearity between the covariates. 

Moreover, it was observed that a reasonable correlation between covariates has little effect on the selection performance of the proposed procedure. However, when the level of correlation becomes high, the performance decreases. This could be improved if structural information on the covariates were \textit{a priori} known. Indeed, in this article, the i.i.d. Bernoulli prior on the indicators $\boldsymbol{\delta}$ \eqref{prior_delta}, entails the assumption that each covariate has the same probability a priori of being included in the model. 
However, there are situations, such as genomic data, in which certain covariates are \textit{a priori} more likely to be included together in the model. This \textit{a priori} structural information on the covariates can be taken into account in SAEMVS by choosing a more flexible prior on $\boldsymbol{\delta}$.
In \cite{stingo_bayesian_2010} and \cite{stingo2011variable}, authors propose the independent logistic regression prior or the Markov random field prior. This could also be considered in our methodology.

Another important remark is that, in this article, we considered a Gaussian distribution for $p(\boldsymbol{y}\vert\boldsymbol{\varphi},\sigma^2)$ in model~\eqref{model}. It is possible to relax this assumption and consider larger distribution classes, such as discrete distributions like the Poisson distribution for example. The proposed methodology can therefore be applied in many contexts.

\section*{Declarations}

\subsubsection*{Funding}
This work was funded by the Stat4Plant project ANR-20-CE45-0012. 

\subsubsection*{Conflict of interest}
The authors declare that they have no conflict of interest.

\subsubsection*{Data and code availability}
Code for running the simulations and generating exact copies of all the figures in this paper, and the real data-set are available at: \\
\href{https://github.com/Marion-Naveau/Supp_Information_SAEMVS}{https://github.com/Marion-Naveau/Supp\_Information\_SAEMVS}

\bibliography{sn-bibliography}


\newpage
\begin{appendices}

\renewcommand\thefigure{\arabic{figure}}
\setcounter{figure}{6}

\section{Algorithms: synthesis and implementation details}
\subsection{MCMC-SAEM algorithm in SSNLME model}
\label{app_MCMCSAEM}

First, in Subsection~\ref{Decomposition}, we notice that $\overset{\sim}{Q}_{1}(\boldsymbol{y},\boldsymbol{\varphi},\theta,\Theta^{(k)})$ takes an exponential form. More precisely, we have:
\begin{equation}
\label{expo}
    \overset{\sim}{Q}_{1}(\boldsymbol{y},\boldsymbol{\varphi},\theta,\Theta^{(k)})= - \psi(\theta,\Theta^{(k)}) + \bigg\langle S(\boldsymbol{y},\boldsymbol{\varphi}), \phi(\theta) \bigg\rangle,
\end{equation}
with:
\begin{itemize}
    \item $S(\boldsymbol{y},\boldsymbol{\varphi}) = \left(\sum_{i,j} (y_{ij}-g(\varphi_i,t_{ij}))^2 \text{ , } vec(\boldsymbol{\varphi}^\top\boldsymbol{\varphi}) \text{ , } vec(\boldsymbol{\varphi})\right)$
    \item $\phi(\theta)=\left(-\dfrac{1}{2\sigma^2} \text{ , } -\dfrac{1}{2}vec(\Gamma^{-1}) \text{ , } vec(\Tilde{V}\boldsymbol{\Tilde{\beta}}\Gamma^{-1})\right)$
    \item $\psi(\theta,\Theta^{(k)})= \dfrac{1}{2} \langle (\Tilde{V}\boldsymbol{\Tilde{\beta}})^\top\Tilde{V}\boldsymbol{\Tilde{\beta}}, \Gamma^{-1} \rangle +\dfrac{1}{2}\sum_{m=1}^q\sum_{\ell'=1}^{p+1} \Tilde{\beta}_{\ell' m}^2 \overset{\sim}{d}_{\ell' m}^{*}(\Theta^{(k)})+\dfrac{n_{tot} + \nu_{\sigma}+2}{2}\log(\sigma^2) + \dfrac{n+d+q+1}{2}\log(\vert\Gamma\vert) + \dfrac{\nu_{\sigma} \lambda_{\sigma}}{2\sigma^2} + \dfrac{1}{2}\mathrm{Tr}(\Sigma_{\Gamma} \Gamma^{-1}) $
\end{itemize}
where $vec(A)$ denotes the vectorisation of a matrix $A$. To simplify the formulas and using that for two matrix $A$ and $B$, $\langle A, B \rangle := \mathrm{Tr}(A^\top B) = \langle vec(A), vec(B) \rangle$, we denote by 
\begin{equation}
\label{S_matriciel}
    (s_1(\boldsymbol{y},\boldsymbol{\varphi}),s_2(\boldsymbol{\varphi}),s_3(\boldsymbol{\varphi}))=\left(\sum_{i,j} (y_{ij}-g(\varphi_i,t_{ij}))^2 \text{ , } \boldsymbol{\varphi}^\top\boldsymbol{\varphi} \text{ , } \boldsymbol{\varphi}\right).
\end{equation}

The use of the decomposition discussed in Subsection~\ref{Decomposition} leads to the following extension of the MCMC-SAEM algorithm, Algorithm~\ref{algoMCMCSAEM}, for computing the MAP estimator of $\Theta$ in the SSNLME model~\eqref{model}~-~\eqref{priors}, where $\Lambda_{\theta}$ denotes the parameter space restricted to $\theta$, $h$ is small (between $1$ and $5$), and $K$ is usually in the order of a few hundred.

In practice, to allow more flexibility during the first iterations and thus to move away more quickly from the initial condition, it is usual to start the algorithm with $n_\textrm{burnin}$ burn-in iterations, \textit{i.e.} to use a step sizes sequence $(\gamma_k)_k$ of the form: $\gamma_k=1$ for $0\leq k \leq n_\textrm{burnin}-1$ and $\gamma_k=1/(k-n_\textrm{burnin}+1)^{\gamma}$ for $n_\textrm{burnin}\leq k \leq K-1$, where $\gamma \in ]0.5,1[$, $n_\textrm{burnin} < K $ with $K$ the number of iterations of the SAEM algorithm \citep[see][]{kuhn2005}. Moreover, to avoid convergence toward a local maximum in the MCMC-SAEM algorithm, a simulated annealing version of SAEM \citep[see][]{lavielle} is implemented.

\begin{algorithm}[htbp]
\caption{MCMC-SAEM}\label{algoMCMCSAEM}
{\small
\begin{algorithmic}
 \State {\bfseries Input:} $K \in \mathbb{N}^*$, $\Theta^{(0)}$ initial parameter, hyperparameters vector $\Xi~=~(\nu_0, \nu_1, \sigma_{\mu}^2, \nu_{\sigma}, \lambda_{\sigma}, \Sigma_{\Gamma}, d, a, b)$, $S_0=0$ and $(\gamma_k)_k$ a step sizes sequence decreasing towards 0 such that $\forall k$, $\gamma_k\in [0,1]$, $\sum_k \gamma_k=\infty$ and $\sum_k \gamma_k^2<\infty$. \\
\For{$k = 0$ to $K-1$} 
\State \begin{enumerate}
    \item \textbf{S-Step:} simulate $\boldsymbol{\varphi}^{(k)}$ using the result of $h$ iterations of an MCMC procedure with $\pi(\boldsymbol{\varphi}\vert\boldsymbol{y},\Theta^{(k)})$ for target distribution.
    \item \textbf{SA-Step:} for $u \in \{1,2,3\}$, compute $s_{u,k+1}=s_{u,k} + \gamma_k(s_u(\boldsymbol{y},\boldsymbol{\varphi}^{(k)})-s_{u,k})$ with $s_u(\boldsymbol{y},\boldsymbol{\varphi}^{(k)})$ defined by \eqref{S_matriciel}.
    \item \textbf{M-Step:} update $\theta^{(k+1)}= \underset{\theta \in \Lambda_{\theta}}{\operatorname{argmax}} \text{ } \Bigg\{ -\psi(\theta,\Theta^{(k)}) + \langle S_{k+1}, \phi(\theta) \rangle \Bigg\}$ and $\alpha^{(k+1)}~=~\underset{\alpha \in [0,1]^q}{\operatorname{argmax}} \text{ }\overset{\sim}{Q}_{2}(\alpha,\Theta^{(k)})$, which reduces to the following explicit forms in model~\eqref{model}-\eqref{priors}:
    \begin{itemize}
    \item $vec(\boldsymbol{\Tilde{\beta}}^{(k+1)})=\left( I_q \otimes \Tilde{V}^\top\Tilde{V} + (\Gamma^{(k)} \otimes I_{p+1})\text{diag}\left(vec(\Tilde{D}^{*})\right)\right)^{-1} vec( \Tilde{V}^\top s_{3,k+1})$
    \item $\Gamma^{{(k+1)}}=\dfrac{\Sigma_{\Gamma}+s_{2,k+1}- (\Tilde{V}\boldsymbol{\Tilde{\beta}}^{(k+1)})^\top s_{3,k+1} - s_{3,k+1}^\top \Tilde{V}\boldsymbol{\Tilde{\beta}}^{(k+1)}  +(\Tilde{V}\boldsymbol{\Tilde{\beta}}^{(k+1)})^\top\Tilde{V}\boldsymbol{\Tilde{\beta}}^{(k+1)}}{n+d+q+1}$
    \item $\sigma^{2^{(k+1)}}=\dfrac{\nu_{\sigma}\lambda_{\sigma}+s_{1,k+1}}{n_{tot}+\nu_{\sigma}+2}$ 
    \item $\alpha_m^{(k+1)}=\dfrac{\sum_{\ell=1}^{p} p_{\ell m}^{*}(\Theta^{(k)})+a_m-1}{p+b_m+a_m-2}$ for $1 \leq m \leq q$
    \end{itemize}
     where $\psi$ and $\phi$ are defined by \eqref{expo}, and $\overset{\sim}{Q}_{2}(\alpha,\Theta^{(k)})$, $(p_{\ell m}^{*}(\Theta^{(k)}))_{1 \leq \ell \leq p ; 1 \leq m \leq q}$ and $(\overset{\sim}{d}_{\ell' m}^{*}(\Theta^{(k)}))_{1 \leq \ell' \leq p+1 ; 1 \leq m \leq q}$ are defined in Proposition~\ref{PropQtilde}, with $\Tilde{D}^{*}~=~(\overset{\sim}{d}_{\ell' m}^{*}(\Theta^{(k)}))_{1\leq \ell' \leq p+1 ; 1 \leq m \leq q} \in \mathcal{M}_{(p+1) \times q}$.
\end{enumerate}
\EndFor
\State {\bfseries Output:} $\widehat{\Theta}^{MAP}=(\widehat{\boldsymbol{\Tilde{\beta}}}^{MAP},\widehat{\Gamma}^{MAP},\widehat{\sigma}^{2^{MAP}}, \widehat{\alpha}^{MAP})=(\boldsymbol{\Tilde{\beta}}^{(K)},\Gamma^{(K)},\sigma^{2^{(K)}},\alpha^{(K)})$.
\end{algorithmic}
}
\end{algorithm}

\newpage
\subsection{Proposed variable selection procedure: SAEMVS}
\label{app_SAEMVS}

The proposed variable selection procedure SAEMVS can be summarised as in Algorithm~\ref{algoSAEMVS}.
\vspace{0.1cm}
\begin{algorithm}[htbp]
\caption{SAEMVS procedure}\label{algoSAEMVS}
\begin{algorithmic}
\State {\bfseries Input:} $\Delta$ a grid of $\nu_0$ values, and all required arguments for MCMC-SAEM (Algorithm~\ref{algoMCMCSAEM}). \\
\State {$\triangleright$ \bfseries Reduce the model collection:}
    \For{$\nu_0 \in \Delta$} \begin{enumerate}
        \item Compute the MAP estimate $\widehat{\Theta}^{MAP}_{\nu_0}$ by Algorithm~\ref{algoMCMCSAEM}. 
        \item Threshold the estimator $\widehat{\boldsymbol{\beta}}^{MAP}_{\nu_0}$ to define sub-model $\widehat{S}_{\nu_0}$ according to Equation~\eqref{S_m}.
        \end{enumerate}
    \EndFor \\
\State {$\triangleright$ \bfseries Compute the eBIC criterion:}
    \For{each unique sub-model among $(\widehat{S}_{\nu_0})_{\nu_0 \in \Delta}$} \begin{enumerate}
        \item Compute the MLE estimate $\widehat{\theta}^{MLE}_{\nu_0}$ in sub-model $\widehat{S}_{\nu_0}$ with an MCMC-SAEM algorithm.               
        \item Compute the log-likelihood $\log p(\boldsymbol{y};\hat\theta^{MLE}_{\nu_0})$ with importance sampling techniques.
        \item Compute the associated $\text{eBIC}(\widehat{S}_{\nu_0})$ according to Equation~\eqref{eBIC}.
    \end{enumerate}
    \EndFor \\
\State {$\triangleright$ \bfseries Identify the best level of sparsity:} compute $\hat \nu_0$ defined by Equation~\eqref{nu0_chap}. \\
\State \bfseries Output:
$\widehat{S}_{\hat\nu_0}$.
 
\end{algorithmic}
\end{algorithm}

\newpage
\subsection{Algorithmic settings of SAEMVS for the comparison study}
\label{setting_SAEMVS_multi}
For the comparison study, the following settings are used for Algorithm~\ref{algoSAEMVS}.
\begin{itemize}
\item The hyperparameter values are set to $\nu_{\sigma}=\lambda_{\sigma}=1$, $d=4$, $\Sigma_{\Gamma}=0.2 I_2$, $a=(1,1)^\top$, $b=(p,p)^\top$, $\sigma_{\mu}=5$, $\nu_1=1000$, and the spike parameter $\nu_0$ runs through a grid $\Delta$ defined as $\log_{10}(\Delta)=\bigg\{-3 + k\times \dfrac{1}{3}, k \in \{0, \dots, 9\}\bigg\}$.
\item The step sizes are defined with $\gamma=2/3$, $n_\textrm{burnin}=150$ and $K=300$ as explained in Appendix \ref{app_MCMCSAEM}.
\item The MCMC-SAEM algorithm is initialised with: $\forall m \in \{1,2\}, \forall \ell \in \{1,\dots, 10\}$ $\boldsymbol{\beta}^{(0)}_{\ell m}=1$, and $\forall m \in \{1,2\}, \forall \ell \in \{11,\dots, p\}$ $\boldsymbol{\beta}^{(0)}_{\ell m}=0.1$, $\mu^{(0)}=(10,10)^\top$, $\sigma^{2^{(0)}}=10^{-2}$, $\Gamma^{(0)}~=~\left (
   \begin{array}{cc}
      0.5 & 0.1 \\
      0.1 & 0.5 \\ 
   \end{array}
\right)$ and $\alpha^{(0)}=0.5$. Note that different initialisations have been tested and have shown similar performances. 
\end{itemize}

\begin{figure}[ht]
  \begin{center}
    \subfloat{
      \includegraphics[width=0.5\textwidth]{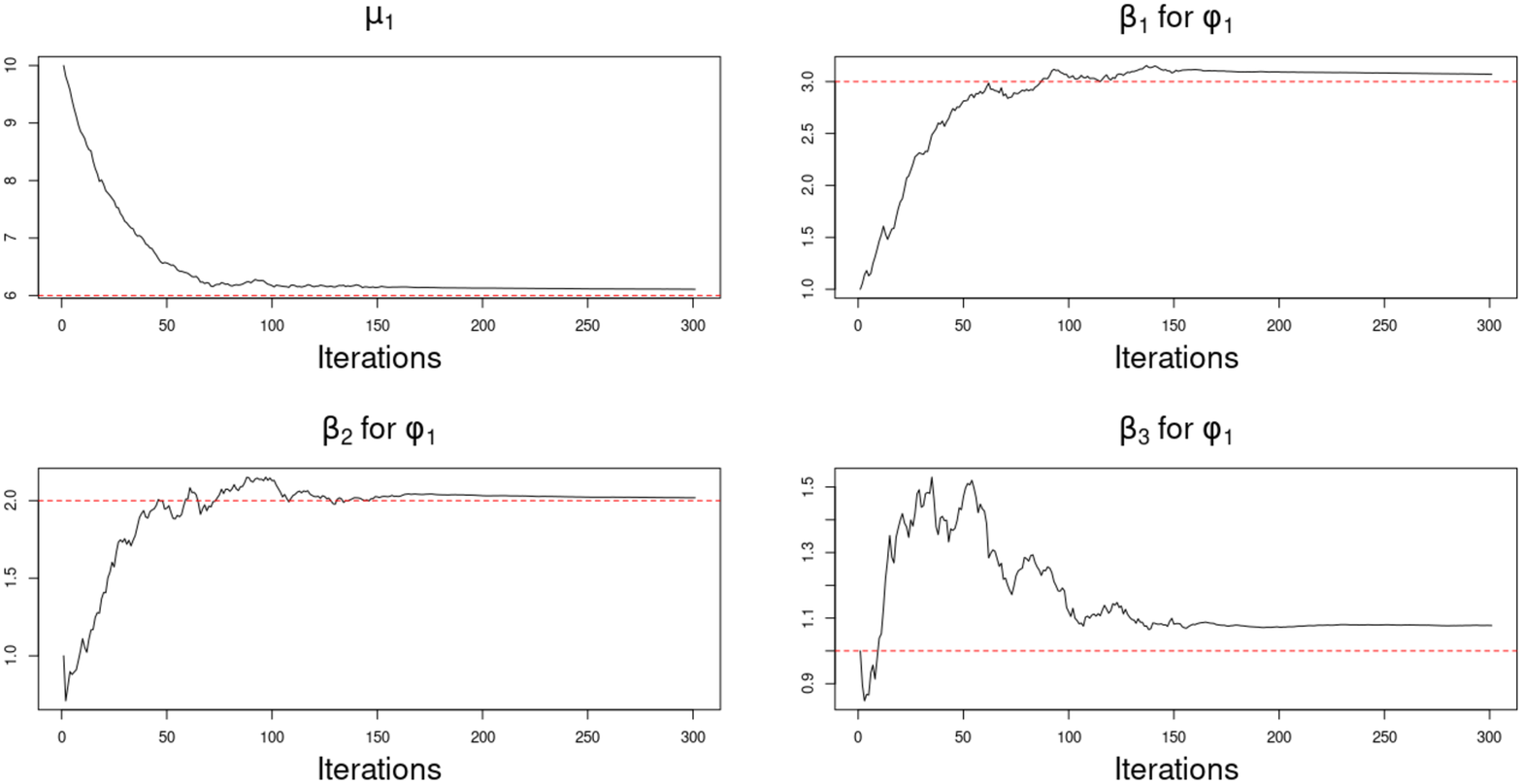}
            }
    \subfloat{
      \includegraphics[width=0.5\textwidth]{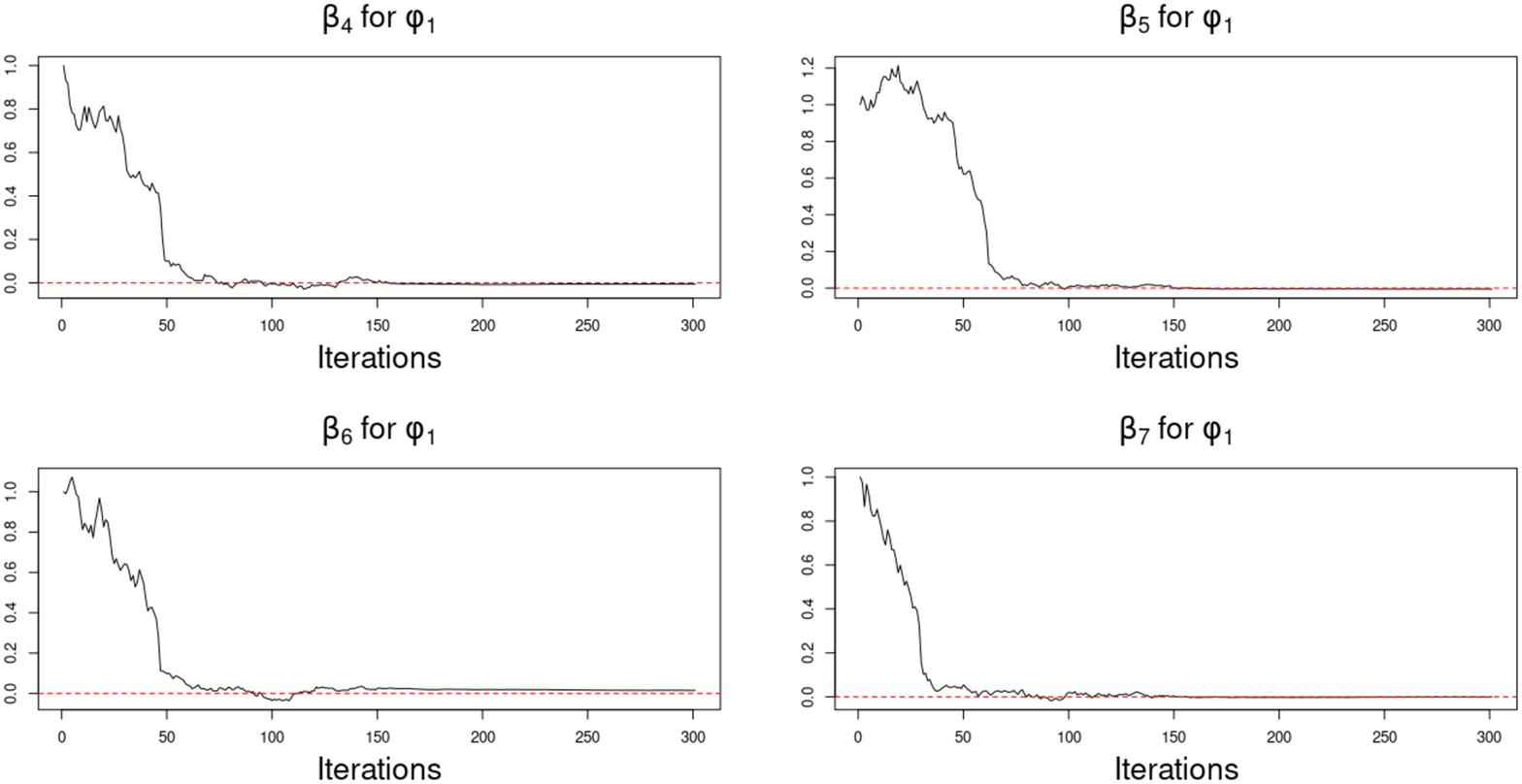}
             } \\
    \subfloat{
      \includegraphics[width=0.5\textwidth]{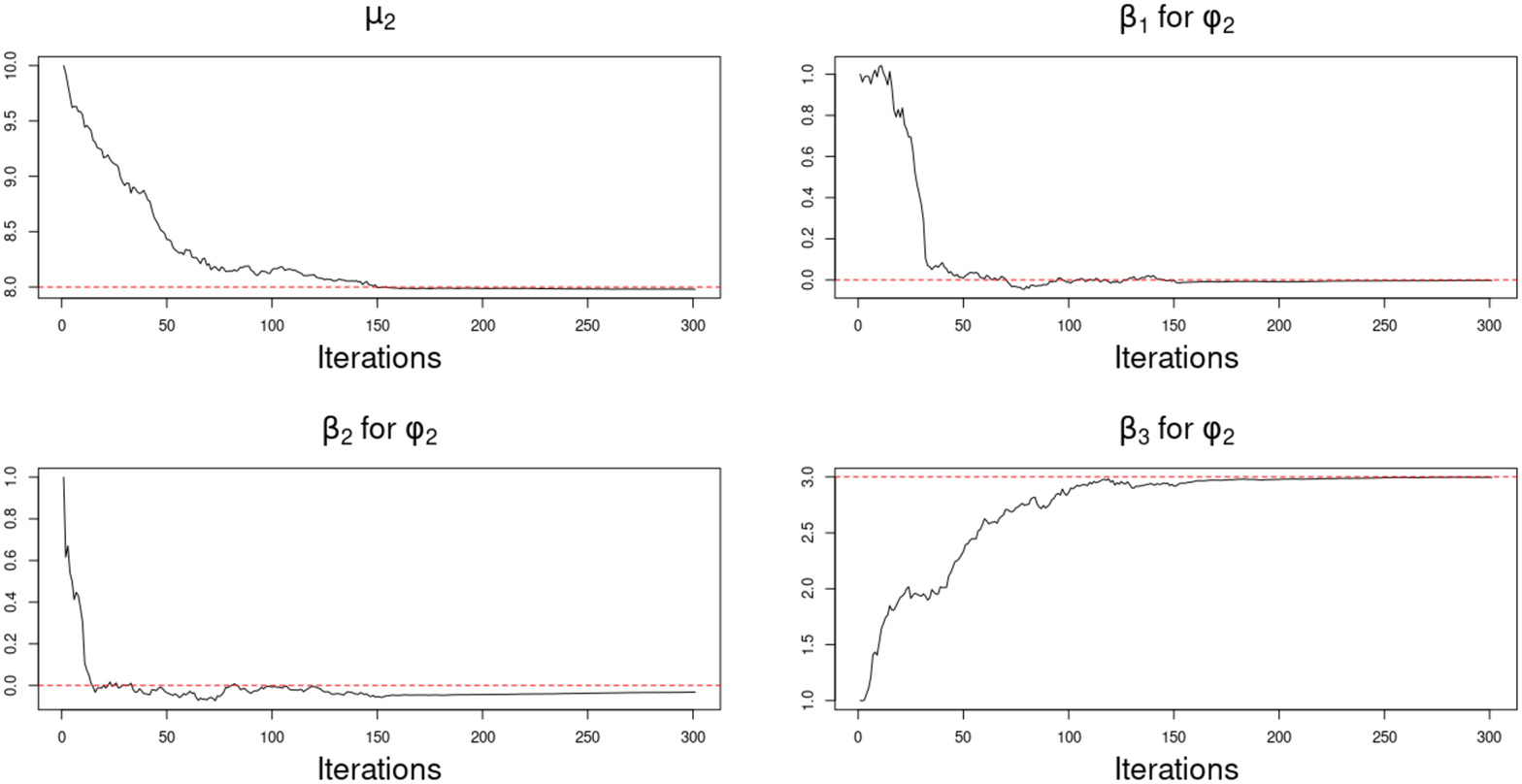}
             }
    \subfloat{
      \includegraphics[width=0.5\textwidth]{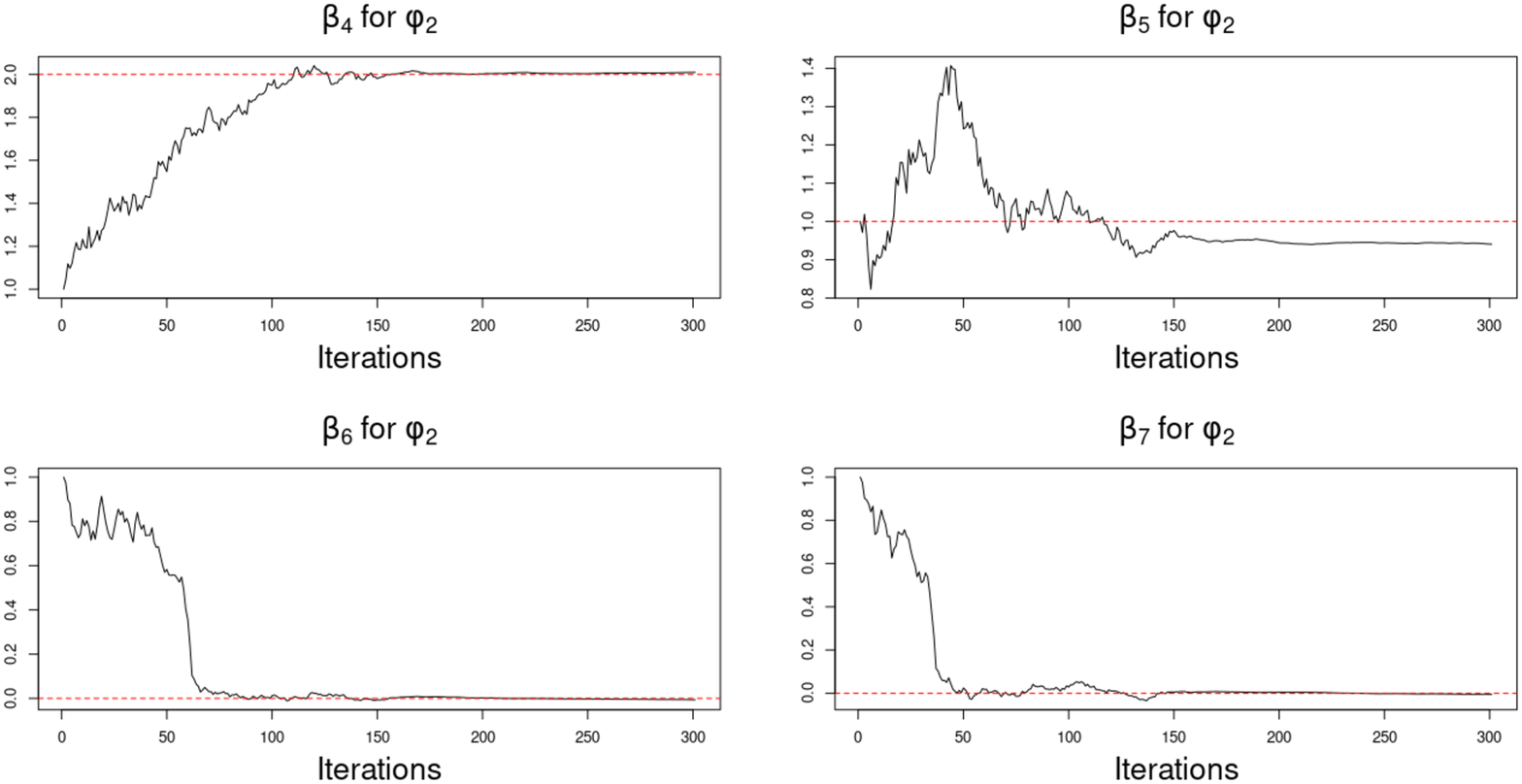}
             } \\
    \subfloat{
      \includegraphics[width=0.5\textwidth]{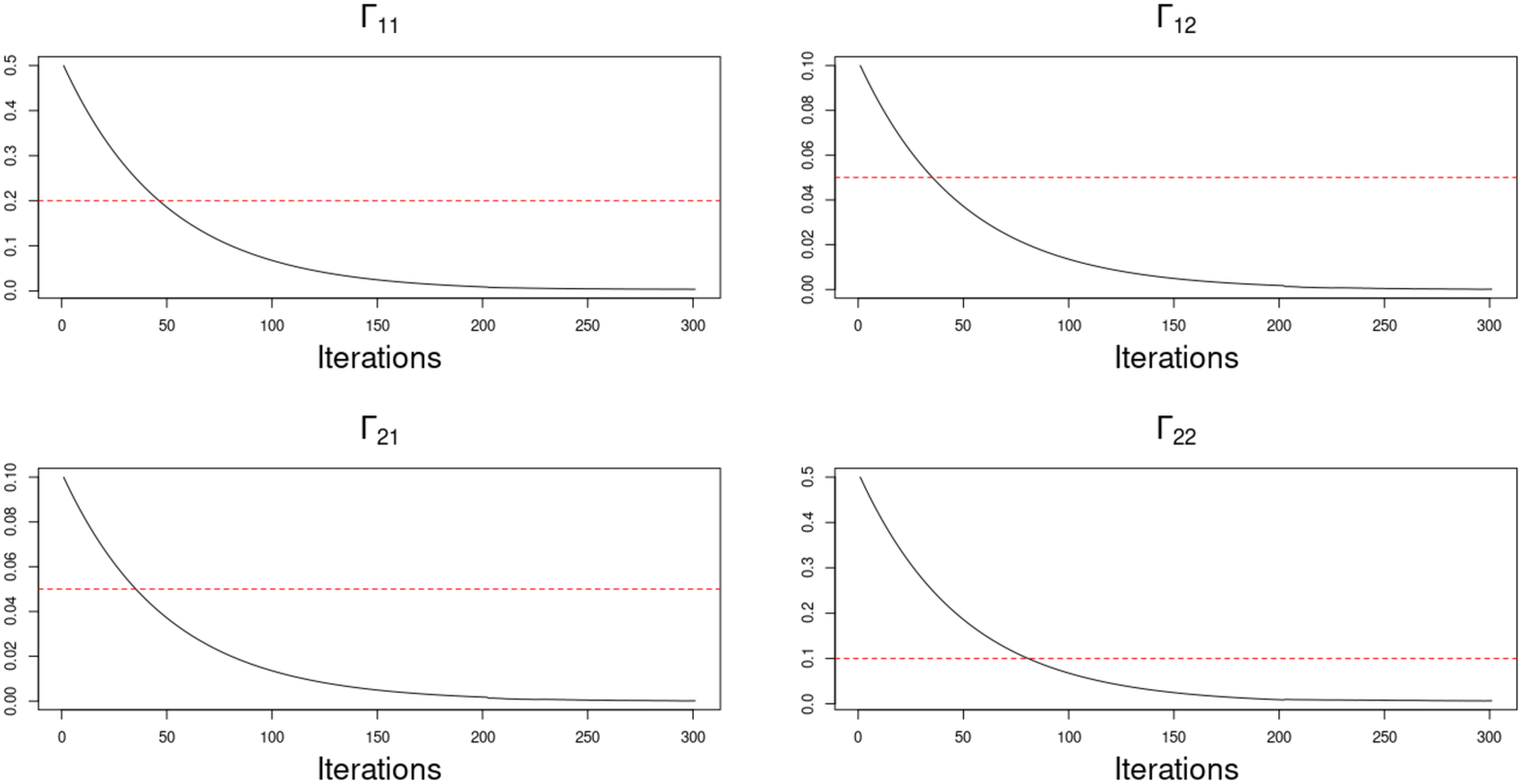}
            }
    \subfloat{
      \includegraphics[width=0.5\textwidth]{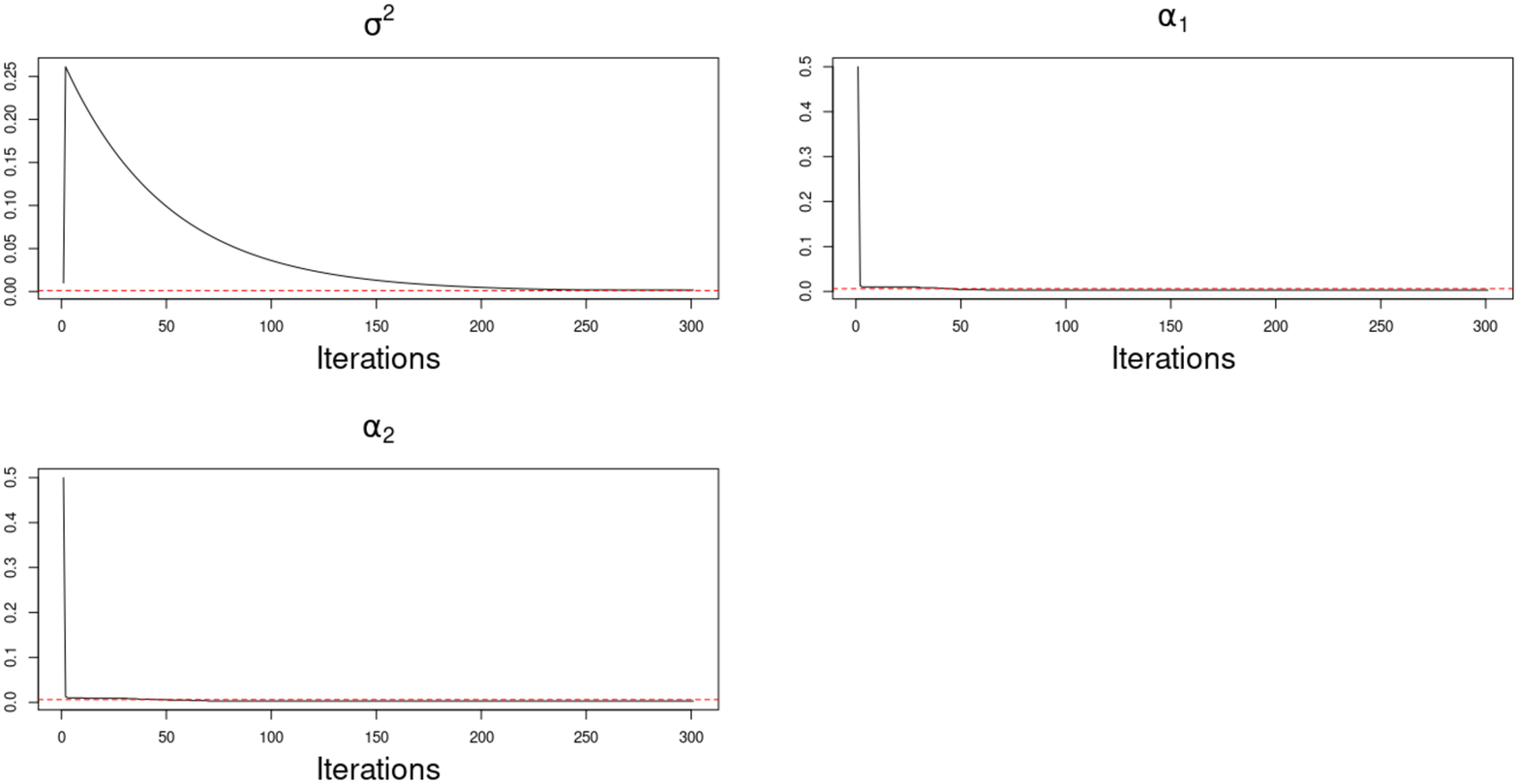}
             } 
    \caption{Convergence graphs of the MCMC-SAEM algorithm for $\mu$, some components of $\boldsymbol{\beta}$, $\sigma^2$, $\Gamma$ and $\alpha$ on one simulated data-set, for $\nu_0=0.005$ and $\nu_1=1000$. The red dashed line corresponds to the true value of the considered parameter.}
    \label{graphcv}
  \end{center}
\end{figure}

Figure \ref{graphcv} represents the convergence graphs of one run of the MCMC-SAEM algorithm for $\mu$, some components of $\boldsymbol{\beta}$, $\sigma^2$, $\Gamma$ and $\alpha$. It is observed that the algorithm converges in a few iterations for any parameter. Note that the parameters are all relatively correctly estimated, except for $\Gamma$ but this was expected because of a over-fitting situation. Indeed, the underestimation of $\Gamma$ can be explained by the fact that since $\nu_0>0$, none of the estimates of the coefficients of $\boldsymbol{\beta}$ is zero and therefore all the covariates are active in the model, which makes the variance estimation of the random effect tend towards $0$.

\subsection{Extension in one-dimensional setting with estimation of fixed effects}
\label{app_extension_fixed}

As it is explained in Subsection~\ref{sim_model}, SAEMVS can easily be adapted to a model where fixed effects must be estimated. For this, the following general one-dimensional model is considered:
\begin{equation*}
    \left\{
    \begin{array}{ll}
        y_{ij} & \overset{\text{ind.}}{\sim} \mathcal{N}(g(\varphi_i, \psi, t_{ij}),\sigma^2), \\
        \varphi_i & \overset{\text{ind.}}{\sim} \mathcal{N}(\mu + \beta^\top V_i,\Gamma^2),
    \end{array}
\right.
\end{equation*}
where $\psi \in \mathbb{R}^s$ are fixed effects to be estimated. Thus, by considering $\psi$ as a latent variable with a normal distribution centred in $\eta$, an unknown parameter, and with a known covariance matrix $\Omega$, and no longer as a parameter, it is possible to obtain an exponential form similar to Equation~\eqref{expo} for $\overset{\sim}{Q}_{1}$. The priors described in Subsection~\ref{sim_model} in the particular case of the logistic growth model are used here, and we have $\Theta~=~(\mu,\beta,\eta,\sigma^2,\Gamma^2,\alpha)$ and $Z=(\varphi,\psi,\delta)$. The calculation of the quantity $Q$ of the EM algorithm becomes: 

\begin{align*}
    Q(\Theta\vert\Theta^{(k)})&= \mathbb{E}_{(\varphi,\psi,\delta) \vert(\boldsymbol{y},\Theta^{(k)})}[\log(\pi(\Theta,\varphi, \psi, \delta\vert\boldsymbol{y}))\vert\boldsymbol{y},\Theta^{(k)}] \\
    &= \mathbb{E}_{(\varphi,\psi)\vert(\boldsymbol{y},\Theta^{(k)})}\left[\overset{\sim}{Q}(\boldsymbol{y},\varphi,\psi,\Theta,\Theta^{(k)})\bigg\vert\boldsymbol{y},\Theta^{(k)}\right],
\end{align*}
where:
\begin{align*}
    \overset{\sim}{Q}(\boldsymbol{y},\varphi,\psi,\Theta,\Theta^{(k)})&= \mathbb{E}_{\delta\vert(\varphi,\psi, \boldsymbol{y},\Theta^{(k)})}[\log(\pi(\Theta,\varphi,\psi,\delta\vert\boldsymbol{y})) \vert \varphi,\psi, \boldsymbol{y},\Theta^{(k)}] \\
    &=  C + \overset{\sim}{Q}_{1}(\boldsymbol{y},\varphi,\psi,\theta,\Theta^{(k)}) + \overset{\sim}{Q}_{2}(\alpha,\Theta^{(k)}),
\end{align*}
with
\begin{align*}
    \overset{\sim}{Q}_{1}(\boldsymbol{y},\varphi,\psi,\theta,\Theta^{(k)})=&-\dfrac{1}{2\sigma^2} \sum_{i,j} (y_{ij}-g(\varphi_i,\psi,t_{ij}))^2 -\dfrac{1} {2\Gamma^2} \vert\vert\varphi - \Tilde{V}\Tilde{\beta}\vert\vert^2-\\
    &\dfrac{1}{2} \sum_{\ell'=1}^{p+1} \Tilde{\beta}_{\ell'}^2 \overset{\sim}{d_{\ell'}^{*}}(\Theta^{(k)})- 
    \dfrac{n_{tot}+\nu_{\sigma}+2}{2}\log(\sigma^2)-\dfrac{n+\nu_{\Gamma}+2}{2}\log(\Gamma^2) -\\ &\dfrac{\nu_{\Gamma}\lambda_{\Gamma}}{2\Gamma^2} -
    \dfrac{\nu_{\sigma}\lambda_{\sigma}}{2\sigma^2}  - \sum_{r=1}^s \dfrac{(\psi_r-\eta_r)^2}{2 \omega_r^2} - \sum_{r=1}^s \dfrac{\eta_r^2}{2 \rho_r^2}
\end{align*}
and $$\overset{\sim}{Q}_{2}(\alpha,\Theta^{(k)})=\log\left(\sqrt{\dfrac{\nu_0}{\nu_1}} \dfrac{\alpha}{1-\alpha}\right)\sum_{\ell=1}^{p} p_{\ell}^{*}(\Theta^{(k)}) + (a-1)\log(\alpha) + (p+b-1)\log(1-\alpha).$$
$(p_{\ell}^{*}(\Theta^{(k)}))_{1 \leq \ell \leq p}$ and $(\overset{\sim}{d_{\ell'}^{*}}(\Theta^{(k)}))_{1 \leq \ell' \leq p+1}$ are defined in Proposition~\ref{PropQtilde}, Equations \eqref{p_l} and \eqref{d_l} with $q=1$.

Note that $\overset{\sim}{Q}_{1}$ is still of the exponential form. Indeed, 
\begin{equation}
\label{Q1tilde}
    \overset{\sim}{Q}_{1}(\boldsymbol{y},\varphi,\psi,\theta,\Theta^{(k)})= - \Psi(\theta,\Theta^{(k)}) + \bigg\langle S(\boldsymbol{y},\varphi,\psi), \phi(\theta) \bigg\rangle
\end{equation}
with:
\begin{itemize}
    \item $S(\boldsymbol{y},\varphi,\psi)=\left(\sum_{i,j} (y_{ij}-g(\varphi_i,\psi,t_{ij}))^2 \text{ , } \sum_{i=1}^{n}\varphi_i^2 \text{ , } \varphi \text{ , } \psi^2 \text{ , } \psi \right)$ 
    \item $\phi(\theta)=\left(-\dfrac{1}{2\sigma^2} \text{ , } -\dfrac{1}{2\Gamma^2} \text{ , } \dfrac{\Tilde{V}\Tilde{\beta}}{\Gamma^2}, \left( -\dfrac{1}{2 \omega_r^2}\right)_{1\leq r\leq s},\left( \dfrac{\eta_r}{\omega_r^2}\right)_{1 \leq r\leq s} \right)$
    \item $\Psi(\theta,\Theta^{(k)})= \dfrac{\vert\vert\Tilde{V}\Tilde{\beta}\vert\vert^2}{2\Gamma^2} +\dfrac{1}{2}\sum_{\ell'=1}^{p+1} \Tilde{\beta}_{\ell'}^2 \overset{\sim}{d_{\ell'}^{*}}(\Theta^{(k)})+\dfrac{n_{tot} + \nu_{\sigma}+2}{2}\log(\sigma^2) + \\ \dfrac{n+\nu_{\Gamma}+2}{2}\log(\Gamma^2) + \dfrac{\nu_{\Gamma} \lambda_{\Gamma}}{2\Gamma^2} + \dfrac{\nu_{\sigma} \lambda_{\sigma}}{2\sigma^2} +\sum_{r=1}^s\dfrac{\eta_r^2}{2\omega_r^2}+\sum_{r=1}^s \dfrac{\eta_r^2}{2 \rho_r^2}$
\end{itemize}

The $k$-th iteration of the MCMC-SAEM algorithm on this model is therefore:

\begin{enumerate}
    \item \textbf{S-Step:} simulate ($\varphi^{(k)}$,$\psi^{(k)}$) using the result of some iterations of a Metropolis-Hastings within Gibbs algorithm with $\pi(\varphi,\psi\vert\boldsymbol{y},\Theta^{(k)})$ for target distribution.
    \item \textbf{SA-Step:} compute $S_{k+1}=S_k + \gamma_k(S(\boldsymbol{y},\varphi^{(k)},\psi^{(k)})-S_k)$ with $S(\boldsymbol{y},\varphi,\psi)$ defined by \eqref{Q1tilde}, where $S_k~=~(s_{1,k},s_{2,k},s_{3,k},s_{4,k},s_{5,k})~\in~ \mathbb{R}\times\mathbb{R}\times\mathbb{R}^n\times\mathbb{R}^q\times\mathbb{R}^q$.
    \item \textbf{M-Step:} update $\theta^{(k+1)}= \underset{\theta \in \Lambda_{\theta}}{\operatorname{argmax}} \text{ } \Bigg\{ -\psi(\theta,\Theta^{(k)}) + \langle S_{k+1}, \phi(\theta) \rangle \Bigg\}$ and \\$\alpha^{(k+1)}~=~\underset{\alpha \in [0,1]}{\operatorname{argmax}} \text{ }\overset{\sim}{Q}_{2}(\alpha,\Theta^{(k)})$. More precisely, \begin{itemize}
        \item ${\Tilde{\beta}^{(k+1)}}=(\Tilde{V}^\top\Tilde{V} +\Gamma^{2^{(k)}}\text{diag}((\overset{\sim}{d_{\ell'}^{*}}(\Theta^{(k)}))_{1\leq \ell' \leq p+1}))^{-1}  \Tilde{V}^\top s_{3,k+1}$,
        \item $\Gamma^{2^{(k+1)}}=\dfrac{\vert\vert\Tilde{V}\Tilde{\beta}^{(k+1)}\vert\vert^2+\nu_{\Gamma} \lambda_{\Gamma}+s_{2,k+1}-2\langle  s_{3,k+1},\Tilde{V}\Tilde{\beta}^{(k+1)}\rangle}{n+\nu_{\Gamma}+2}$,
        \item $\sigma^{2^{(k+1)}}=\dfrac{\nu_{\sigma}\lambda_{\sigma}+s_{1,k+1}}{n_{tot}+\nu_{\sigma}+2}$, 
        \item $\eta_r^{(k+1)}=\dfrac{(s_{5,k+1})_r}{1 +\dfrac{\omega_r^{2}}{\rho_r^2}}$ for $1 \leq r \leq s$,
        \item $\alpha^{(k+1)}=\dfrac{\sum_{\ell=1}^{p} p_{\ell}^{*}(\Theta^{(k)})+a-1}{p+b+a-2}$.
    \end{itemize}
\end{enumerate}

As you can see, $\overset{\sim}{Q}_{1}$ is separable into $(\mu,\beta,\sigma^2,\Gamma^2,\alpha)$ and $\eta$, which means that the inference method used for the parameters $(\mu,\beta,\sigma^2,\Gamma^2,\alpha)$ is unchanged, \textit{i.e.} the formulas to update these parameters in M-step are identical, it is only the way to simulate the sufficient statistics that has changed. 

Thus, thanks to this algorithm, it is obtained an estimation $\hat\theta_{\nu_0}^{MAP}$, where \\ $\hat\theta_{\nu_0}^{MAP}=(\widehat{\mu}^{MAP}_{\nu_0},\widehat{\beta}^{MAP}_{\nu_0},\widehat{\eta}^{MAP}_{\nu_0},\widehat{\Gamma}_{\nu_0}^{2,MAP},\widehat{\sigma}_{\nu_0}^{2,MAP})$, and the estimation of $\eta$ is used as an estimation of $\psi$. Then, to finish the model collection reduction step of SAEMVS, Algorithm~\ref{algoSAEMVS}, the estimator $\widehat{\beta}^{MAP}_{\nu_0}$ is thresholded to obtain a promising sub-model $\widehat{S}_{\nu_0}$ given by Equation~\eqref{S_m}. The selection threshold formula is unchanged because it only depends on the second layer of the model~\eqref{model_simul}.

For the model selection step, to compute the eBIC criterion, it is also necessary to go through the extended model~\eqref{extended_model}. Indeed, as described in \cite{kuhn2005}, the MLE in the sub-model $\widehat{S}_{\nu_0}$ is computed in the extended model by using an MCMC-SAEM algorithm and the estimation of $\eta$ is used as an estimation of $\psi$. Then, the log-likelihood is approached by a Monte-Carlo method: for $T$ large enough,  
$$\small {\log \left(p(\boldsymbol{y};\hat\theta^{MLE}_{\nu_0})\right) \approx \sum_{i=1}^{n} \log\left( \dfrac{(2\pi \widehat{\sigma}^{2^{MLE}}_{\nu_0})^{-n_i/2}}{T} \sum_{t=1}^{T} \exp\left( -\sum_{j=1}^{n_i} \dfrac{(y_{ij}-g( \varphi_i^{(t)},\widehat{\psi}_{\nu_0}^{MLE},t_{ij}))^2}{2\widehat{\sigma}^{2^{MLE}}_{\nu_0}}\right) \right)}$$
where $p(\boldsymbol{y} ; \theta)$ denotes the likelihood of model~\eqref{model_simul}, and for all $i\in \{1,\dots n\}$, $(\varphi_i^{(t)})_{t \in \{1,\dots T\}}$ are simulated i.i.d. according to $p(\varphi_i ; \hat\theta^{MLE}_{\nu_0})~=~\mathcal{N}~(~\widehat{\mu}^{MLE}_{\nu_0}~ +~(\widehat{\beta}_{\nu_0}^{MLE})^\top V_i,\widehat{\Gamma}^{2^{MLE}}_{\nu_0})$.

\subsection{Algorithmic settings in the simulation study}
\label{algo_setting}
For the simulation study, the following settings are used for Algorithm~\ref{algoSAEMVS}.
\begin{itemize}
\item The hyperparameter values are set to $\nu_{\sigma}=\lambda_{\sigma}=\nu_{\Gamma}=\lambda_{\Gamma}=1$, $a=1$, $b=p$, $\sigma_{\mu}=3000$, $\rho_1^2=\rho_2^2=1200$, $\nu_1=12000$, and the spike parameter $\nu_0$ runs through a grid $\Delta$ defined as $\log_{10}(\Delta)=\bigg\{-2 + k\times \dfrac{4}{19}, k \in \{0, \dots, 19\}\bigg\}$.
\item The step sizes are defined with $\gamma=2/3$, $n_\textrm{burnin}=350$ and $K=500$ as explained in Appendix \ref{app_MCMCSAEM}.
\item The MCMC-SAEM algorithm is initialised with: $\forall \ell \in \{1,\dots, 10\}$ $\beta^{(0)}_{\ell}=100$, $\forall \ell~\in~\{11,\dots, p\}$ $\beta^{(0)}_{\ell}=1$, $\mu^{(0)}=1400$, $\sigma^{2^{(0)}}=100$, $\Gamma^{2^{(0)}}=5000$, $\alpha^{(0)}=0.5$ and $\eta^{(0)}= (400,400)^\top$. Note that different initialisations have been tested and have shown similar performances. 
\item At the beginning of the algorithm, $\Omega=\textrm{diag}(20,20)$ and it is slowly reduced during the iterations as explained in Remark~\ref{rmk_omega} with $\kappa=40$ and $\tau=0.9$.
\end{itemize}

\section{Results for scenarios with correlated covariates}
\label{ann_corr}

In this appendix, we give the results for scenarios not discussed in Section \ref{sim_results}.

\begin{figure}[ht]
    \centering
    \includegraphics[width=\textwidth]{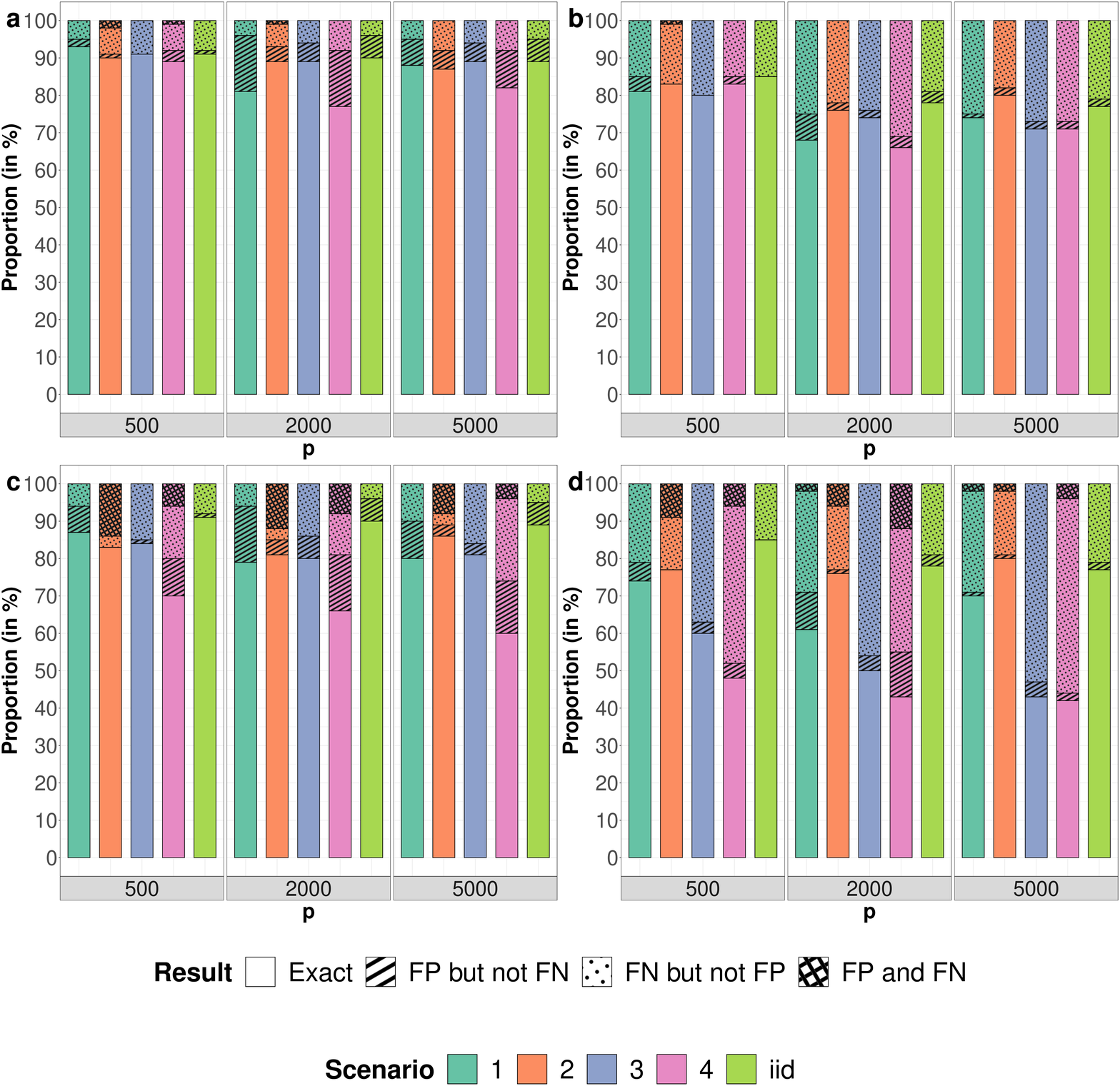}
    \caption{Correlated covariates. Proportion of data-sets on which Algorithm~\ref{algoSAEMVS} selects the correct model ("Exact", unpatterned bars), a model that contains false positives but not false negatives ("FP but not FN", striped bars), a model that contains false negatives but not false positives ("FN but not FP", dotted bars), or a model that contains both false positives and false negatives ("FP and FN", crosshatched bars) for $\rho_{\Sigma}=0.3$ and $\Gamma^2=200$ (a), $\rho_{\Sigma}=0.3$ and $\Gamma^2=2000$ (b), $\rho_{\Sigma}=0.6$ and $\Gamma^2=200$ (c), and $\rho_{\Sigma}=0.6$ and $\Gamma^2=2000$ (d), and different values of $p$. Scenario "i.i.d" corresponds to the case where the covariates are not correlated and is used as a reference.}
\end{figure}

In scenario $2$, it is assumed that the third relevant covariate is correlated to the non-active covariates. In this case, similar results to the i.i.d scenario are observed. Indeed, the selection performances of SAEMVS are only slightly affected by this scenario of correlations. This can be explained by the fact that, in this case, among the group of correlated covariates, the method will tend to select only one (or at least a very limited number of covariates among them): the most intense is chosen, \textit{i.e.} the third true covariate. Next, scenario $3$ describes correlations between the relevant covariates. Like the previous scenario, the procedure tends to select few covariates among the correlated covariates since they explain the response variable in a similar way. This also explains the degradation of the results when $\rho_{\Sigma}$ increases.

\section{Further details on real data}
\label{annexe_data_reel}

\subsection{More details on the variety genotyping process}
\label{app_genotyping}

The varieties of the data-set used in this application was genotyped with the TaBW280K high-throughput genotyping array \citep{rimbert}. This array was designed to cover both genic and intergenic regions of the three bread wheat subgenomes. Markers in strong Linkage Disequilibrium (LD) were filtered out using the pruning function of PLINK \citep{purcell} with a window of size $100$ SNPs (Single Nucleotide Polymorphism, also called "molecular markers" in the following), a step of $5$ SNPs and a LD threshold of $0.8$, as proposed in \cite{charmet}. Missing values were imputed as the marker observed frequency, and then these imputed values are replaced by $0$ or $1$ using a threshold of $0.5$. Monomorphic and unmapped markers were removed from the data-set. Eventually, we obtained $p=$ 26,189 polymorphic high-resolution SNPs with a physical position on the v1 reference genome \citep{international2018shifting}.

\subsection{Representation of the data-set}
\label{app_fig_data_reel}

\begin{figure}[ht]
    \centering
    \includegraphics[width=0.8\textwidth]{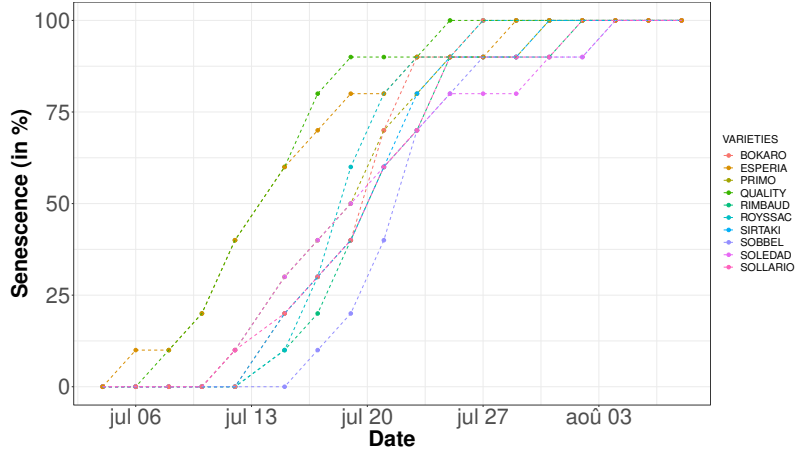}
    \caption{Representation of the data-set for 11 different varieties}
    \label{data_reel}
\end{figure}

Figure~\ref{data_reel} shows part of the real data-set. 

\subsection{Settings}
\label{ann_settings_data_reel}

For application on real data, the following settings are used for Algorithm~\ref{algoSAEMVS}.
\begin{itemize}
\item The hyperparameter values are set to $\nu_{\sigma}=\lambda_{\sigma}=\nu_{\Gamma}=\lambda_{\Gamma}=\nu_{\Omega}=\lambda_{\Omega}=1$, $a=1$, $b=p$, $\sigma^2_{\mu}=\sigma^2_{\lambda}=\sigma^2_{\eta}=100$, $\nu_1=10$, and the spike parameter $\nu_0$ runs through a grid $\Delta$ defined as $\log_{10}(\Delta)=\bigg\{-4.5 + k\times \dfrac{1.5}{9}, k \in \{0, \dots, 9\}\bigg\}$.
\item The step sizes are defined with $\gamma=2/3$, $n_\textrm{burnin}=250$ and $K=400$ as explained in Appendix \ref{app_MCMCSAEM}.
\item The MCMC-SAEM algorithm is initialised with: $(\beta^{(0)})_{\ell}=0.25$ if $\ell$ is a marker less than 1 mega base distance from a heading QTL or a major flowering gene, and otherwise $(\beta^{(0)})_{\ell}=0.1$, $\mu^{(0)}=20$, $\lambda^{(0)}=(0.25,\dots,0.25)^\top$, $\sigma^{2^{(0)}}=80$, $\Gamma^{2^{(0)}}=50$, $\eta^{(0)}=5$, $\Omega^{2^{(0)}}=50$, and $\alpha^{(0)}=0.5$.
\end{itemize}

\subsection{Table of results of the application on real data}
\label{annexe_table}

Table~\ref{Tab_chr} summarises for each chromosome, the number of covariates, the number of heading QTLs, the number of major flowering genes present in that chromosome, and the number of SNPs selected by SAEMVS.

\begin{table}[ht]
\caption{Summary of data and number of SNPs selected by SAEMVS for each chromosome}\label{Tab_chr}%
\begin{tabular}{@{}rrrrr@{}}
\toprule
\textbf{Chromosome} & \textbf{p} & \textbf{Number of} & \textbf{Number of} & \textbf{Number of} \\
 & & \textbf{heading QTLs} & \textbf{flowering genes} & \textbf{selected SNPs} \\
\midrule
1A & 1473 & 1 & 0 & 2\\
1B & 1604 & 0 & 0 & 0\\
1D & 497 & 0 & 0 & 1\\
\midrule
2A & 1416 & 0 & 1 & 3\\
2B & 1672 & 1 & 1 & 2\\
2D & 696 & 7 & 1 & 2\\
\midrule
3A & 1477 & 0 & 0 & 0\\
3B & 1888 & 0 & 0 & 0\\
3D & 722 & 0 & 0 & 1\\
\midrule
4A & 1259 & 0 & 0 & 1\\
4B & 961 & 0 & 0 & 0\\
4D & 571 & 0 & 0 & 0\\
\midrule
5A & 1598 & 0 & 1 & 2\\
5B & 1535 & 0 & 1 & 0\\
5D & 772 & 0 & 1 & 1\\
\midrule
6A & 1119 & 3 & 0 & 1\\
6B & 1317 & 0 & 0 & 0\\
6D & 584 & 0 & 0 & 1\\
\midrule
7A & 1819 & 0 & 1 & 0\\
7B & 1515 & 1 & 1 & 4\\
7D & 778 & 0 & 1 & 1\\
\toprule
\end{tabular}
\end{table}

\end{appendices}


\end{document}